\documentclass[preprint,10pt]{elsarticle}

\usepackage{amssymb}
\usepackage{amscd}
\usepackage{mathrsfs}
\usepackage{amsmath}
\usepackage{amsthm}
\usepackage{mathptmx}
\usepackage{mathtools}
\usepackage{comment}
\usepackage{booktabs}       % professional-quality tables

\usepackage{pifont}

\newtheorem{lemma}{Lemma}[section]
\newtheorem{theorem}{Theorem}[section]

\usepackage{lineno,hyperref}
\modulolinenumbers[5]

\usepackage{soul,color, xcolor}
\soulregister\cite7
\soulregister\ref7

\usepackage{float}
\usepackage{graphicx}
\usepackage{graphics}
\usepackage{subfigure}

\usepackage[a4paper,left=30mm,width=150mm,top=25mm,bottom=25mm]{geometry}

\usepackage{setspace}

\usepackage{indentfirst}
\biboptions{sort&compress}

\usepackage{titlesec} 
\titleformat*{\section}{\large\bfseries}
\titleformat*{\subsection}{\bfseries}
\titleformat*{\subsubsection}{\bfseries}

\begin{document}

\begin{frontmatter}
\title{Implicit-explicit time discretization schemes for a class of semilinear wave equations with nonautonomous dampings\tnoteref{mytitlenote}}
\tnotetext[mytitlenote]{This research was partially supported by the National Natural Science Foundation of China (12272297), the Guangdong Basic and Applied Basic Research Foundation (2022A1515011332) and the Fundamental Research Funds for the Central Universities.}

\author[mymainaddress,secondaddress]{Zhe Jiao\corref{mycorrespondingauthor}}
\cortext[mycorrespondingauthor]{Corresponding author}
%\author[mymainaddress,secondaddress]{Zhe Jiao}
\ead{zjiao@nwpu.edu.cn}

\author[mymainaddress]{Yaxu Li}
\ead{liyaxu919@163.com}

\author[mymainaddress,thirdaddress]{Lijing Zhao}
\ead{zhaolj@nwpu.edu.cn}

\address[mymainaddress]{School of Mathematics and Statistics, Northwestern Polytechnical University, Xi'an 710129, China}

\address[secondaddress]{MOE Key Laboratory for Complexity Science in Aerospace, Northwestern Polytechnical University, Xi'an 710129, China}

\address[thirdaddress]{Research and Development Institute of Northwestern Polytechnical University in Shenzhen, Shenzhen 518057, China }

\begin{abstract}
This paper is concerned about the implicit-explicit (IMEX) methods for a class of dissipative wave systems with time-varying velocity feedbacks and nonlinear potential energies, equipped with different boundary conditions. Firstly, we approximate the problems by using a vanilla IMEX method, which is a second-order scheme for the problems when the damping terms are time-independent. However, rigors analysis shows that  the error rate declines from second to first order due to the nonautonomous dampings. To recover the convergence order, we propose a revised IMEX scheme and apply it to the nonautonomous wave equations with a kinetic boundary condition. Our numerical experiments demonstrate that the revised scheme can not only achieve second-order accuracy but also improve the computational efficiency.
\end{abstract}


\journal{}

\end{frontmatter}
%%%%%%%%%%%%%%%%%%%%%%%%%%%%%%%%%%%%%%%%%%%%%%%%%%%%%%%%%%%%%%%%%%%%%%%%%%%%%%%%

%\linenumbers

\section{Introduction}
\label{sec:intro}

In this paper we consider numerical solutions for the following semilinear evolution equation of second-order
\begin{equation*}
\ddot{u}(t) + B(t)\dot{u}(t) + Au(t) = f(t, u(t))
\end{equation*}
in a suitable Hilbert space. Such evolution equation is renowned for their prevalence in a spectrum of mathematical models derived from physics, elastodynamics, acoustics etc, for example, the propagation of dislocations in crystals~\cite{Kozlowski2006}, wave motions in a dissipative medium~\cite{Chwang1998interaction}, and structural acoustic interactions~\cite{Howe1998acoustics}, thereby providing a versatile abstract framework that encompasses a multitude of wave-type equations wave-type equations with a diverse array of boundary conditions such as Dirichlet, Neumann, Robin and dynamic boundary conditions.

The theoretical exploration of the semilinear evolution equation has attracted many researchers' attention over an extended period. (cf.~\cite{Pucci1993, Pucci1996, liu1998spectral, Cabot2009, vitillaro2013strong, mugnolo2021wave, Jiao2023stability}), while the numerical analysis on this evolution equation can be broadly classified into two main categories based on the boundary conditions. 
%In the first category, it is concerned with the numerical solutions of wave-type equations subject to Dirichlet, Neumann or Robin boundary conditions. There are numerous investigations on this research direction including ...\cite{?}... 
In the first category, it is concerned with wave-type equations subject to Dirichlet, Neumann or Robin boundary conditions. There are numerous investigations on numerical methods for these equations including finite difference methods~\cite{fatkullin2001fdm, deng2020fdm}, finite element methods~\cite{jiang2010fem, he2020fem}, characteristic methods~\cite{bukiet1996characteristic, hasan2022characteristic}, spectral methods~\cite{sternberg1987spectral, fornberg1999fast} and the methods of lines~\cite{saucez2004methodofline, shakeri2008methodofline}, among others.  
The second category focuses on wave equations with dynamic boundary conditions. A series of papers~\cite{Hipp2017, Hipp2019, Hochbruck2020fem, Leibold2021} have been published recently to provide a comprehensive error analysis on some discrete schemes for this type of wave equations. It is noteworthy that the damping terms are time-independent in all these papers. In light of this, this paper is aimed to investigate the numerical schemes for semilinear wave-type equations with time-varying dampings and subject to dynamic boundary conditions, addressing a gap in the current literature.

It is a well-established fact that explicit schemes usually suffer severe time step restrictions despite their computational simplicity, and full implicit schemes are always unconditionally stable but rather taxing and time consuming. Against this backdrop, the focus of this paper is on an IMEX time discretization scheme. This approach is a popular and widely-used approach (cf.~\cite{kadioglu2010imex, lemieux2014imex, akrivis2015imex, gao2016imex}), which leverages an implicit treatment for the unbounded linear components and an explicit treatment for the nonlinear aspects of a differential equation. For‌ an in-depth investigation of the IMEX scheme, we direct the readers to seminal works:~\cite{ascher1997implicit} for IMEX Runge-Kutta schemes, \cite{frank1997imex} for IMEX multistep schemes, \cite{layton2016imex} for Crank-Nicolson-leapfrog IMEX scheme, and references therein.

Closely pertinent to our investigation, we refer to a recent study~\cite{Hochbruck2021imex}, where an IMEX scheme is introduced. This scheme is crafted by merging the Crank-Nicolson method for linear segments and the leapfrog method for nonlinear components, culminating in a second-order numerical scheme. However, when this IMEX scheme is applied to a class of semilinear wave equations with nonautonomous damping terms, it is noted that the error rate descends to the first order, as will be elaborated in Section 4. Thus, it is both natural and necessary to delve into the influence that nonautonomous damping terms exert on the convergence rate of the IMEX scheme.     

Motivated by previous research, this paper aims to present a comprehensive discretization of the aforementioned semilinear evolution equation in an abstract manner. The key contributions of this study are encapsulated as follows:
The major contributions of this article can be summarized as follows: (i) Our methodology offers a cohesive approach to solve a class of wave-type equations with time-dependent dampings and subject to different types of boundary conditions. (ii) We conduct an error analysis on the full discretization with a vanilla IMEX scheme, which reveals the time-varying damping coefficients play a fundamental role in making the error oder decrease to be one. (iii) The full discretization with a revised IMEX scheme is a second-order method that is more computationally efficient without compromising the order of convergence.
We collect the obtained error estimates for the temporal discretization in Table~\ref{comparison}, and mark whether these results are discussed theoretically or illustrated numerically.
\begin{table}[!htbp]
\begin{center}
\begin{tabular}{ccccc}
\toprule
Scheme & autonomous damping & nonautonomous damping & Discussed in & illustrated \\
\midrule
%CN & $\tau^2 $ & $\tau$ & Theorem \ref{main_result} & \ding{51} \\
IMEX & $\tau^2 $ & $\tau$ & Theorem \ref{main_result} & \ding{51}\\
{\color{cyan}Revised IMEX}  & $\tau^2 $ & $\tau^2 $ & \ding{55} & \ding{51} \\
\bottomrule
\end{tabular}
\caption{Overview of convergence rate for time discretization schemes with step size $\tau$.}
\end{center}
\label{comparison}
\end{table}

This paper is structured as follows. In Section 2, we begin by introducing semilinear evolution equations as an abstract framework for wave-type equations with time-dependent damping, and then present three exemplary wave equations with distinct boundary conditions. Section~\ref{sec:full_discretization} commences with an overview of a general space discretization. Subsequently, we develop two full discrete schemes by combing the vanilla and revised IMEX with the abstract space discretization. An error analysis of these full discretizations will be also given. Numerical examples are presented in Section 4 to substantiate the error bounds derived from our main theorem. These examples will not only highlight the error rate but also showcase the efficacy of the proposed schemes. The focus of Section 5 is on the detailed proof of our main result. For better readability, the Crank-Nicolson scheme and the vanilla IMEX scheme for the semilinear wave equations with nonautonomous dampings are detailed in \ref{sec:CN} and \ref{sec:IMEX}, respectively.

%%%%%%%%%%%%%%%%%%%%%%%%%%%%%%%%%%%%%%%%%%%%%%%%%%%%%%%%%%%%%%%%%%%%%%%%%%%%%%%%
\section{Preliminary}
\label{sec:continuous_case}

%In this section, ...

Let $V$, $H$ be Hilbert spaces with norms $\|\cdot\|_V$ and $\|\cdot\|_H$. The second norm induced by an inner product $m(\cdot, \cdot)$ on $H$, and $V$ is densely embedded in $H$. Moreover, we identify $H$ with its dual $H^\ast$ which gives the Gelfand triple 
\[
	V \hookrightarrow H \equiv H^{\ast} \hookrightarrow V^\ast.
\]
Here, the notation $\hookrightarrow$ means the dense embedding. Thus, there exists an embedding constant $C_{\mathrm{em}} > 0$ such that $\|\upsilon\|_H \leqslant C_{\mathrm{em}} \|\upsilon\|_V$, and the duality $\langle \cdot, \cdot\rangle_V: V^\ast\times V\rightarrow\mathbb{R}$ coincide with $m(\cdot, \cdot)$ on $H \times V$. 
%Throughout the paper, $C_1, C_2, \cdots$ denote corresponding constants and $o(\cdot)$ means the higher-order infinitesimal.

The abstract wave equation that covers all examples in this paper can be written as follows.
\begin{equation} \label{eq:secondorder}
\begin{aligned}
	\langle\ddot{u}, \phi\rangle_V + b(\dot{u}, \phi) + \gamma(t) m(\dot{u}, \phi) +a(u, \phi) &= \langle f(t, u), \phi\rangle_V  \quad \textrm{for all $\phi\in V$, and $t\in(0, T]$,} \\
                                        u(0)= u^{0}, \quad \dot{u}(0) &= v^0.
\end{aligned}
\end{equation}
where the initial data $u^0\in V$, $v^0\in H$, $a: V\times V \rightarrow \mathbb{R}$ and $b: V\times V\rightarrow \mathbb{R}$ are continuous bilinear forms, and the nonlinearity $f: [0, T]\times V\rightarrow V^\ast$ is continuous. The notations $\dot{} = \frac{\mathrm{d}}{\mathrm{d}t}$ and $\ddot{} = \frac{\mathrm{d}^2}{\mathrm{d}t^2}$ denote the first and second temporal derivatives, respectively.

We define operators $A, B\in\mathcal{L}(V, V^\ast)$ induced by the bilinear forms $a$ and $b$ via
\begin{equation*} 
\begin{aligned}
            \langle A\varphi, \phi\rangle_V = a(\varphi, \phi) \quad \textrm{and} \quad
            \langle B\varphi, \phi\rangle_V = b(\varphi, \phi)   \quad \textrm{for all $\varphi, \phi\in V$}.
\end{aligned}
\end{equation*}
The variational equation \eqref{eq:secondorder} can be equivalently stated as an second-order evolution equation in $V^\ast$
\begin{equation}\label{eq:wave}
\begin{aligned}
\ddot{u}(t) + B(t)\dot{u}(t) + Au(t) &= f(t, u(t)), \quad t\geqslant 0,\\
u(0) =u^0, \quad \dot{u}(0) &= v^0
\end{aligned}
\end{equation}
with $B(t) = B + \gamma(t)I$.
In addition, if we set $\dot{u} = v$ and define
\begin{equation*}
\mathbf{u}=\left[
\begin{array}{c}
u\\
v
\end{array}
\right],
\quad
S=\left[
\begin{array}{cc}
0 & -I\\
A & B
\end{array}
\right],
\quad
g(t, \mathbf{u})=\left[
\begin{array}{c}
0\\
f(t, u) - \gamma(t)v
\end{array}
\right],
\quad
\mathbf{u}^0=\left[
\begin{array}{c}
u^0\\
v^0
\end{array}
\right],
\end{equation*}
then \eqref{eq:wave} can be reformulated to be the following first-order equation in the Hilbert space $X=V\times H$.
\begin{equation}\label{eq:evolution}
\begin{aligned}
\dot{\mathbf{u}}(t) + S\mathbf{u}(t) &= g(t, \mathbf{u}(t)), \quad t\geqslant 0,\\
\mathbf{u}(0) &= \mathbf{u}^0.
\end{aligned}
\end{equation}

We need the following assumptions to guarantee the evolution equation \eqref{eq:wave} is locally wellposed. 
\begin{itemize}
\item[(H1)] The bilinear form $a(\cdot, \cdot)$ is symmetric and there exists a constant $\alpha>0$ such that
	\[
		\mathrm{p}(\cdot, \cdot) = a(\cdot, \cdot) + \alpha m(\cdot, \cdot)
	\]
	is a inner product on $V$. Let the norm $\|\cdot\|_V$ be generated by this inner product.   
\item[(H2)] The bilinear form $b$ is bounded, that is,
	\[
		b(\upsilon, \omega) \leqslant C_{\mathrm{b}}\|\upsilon\|_V \|\omega\|_H \quad \textrm{for any $\upsilon\in V$, $\omega\in H$}
	\]
	and quasi-monotone which means there exists a constant $\beta>0$ for any $\upsilon\in V$ 
	\[
		b(\upsilon, \upsilon) + \beta m(\upsilon, \upsilon) \geqslant 0.
	\]
\item[(H3)] The nonlinearity $f$ belongs to $C(V; H)$ and also is locally Lipschitz-continuous, that is, for any $\upsilon, \omega \in V$ satisfying $\| \upsilon \|_V, \|\omega \|_V \leqslant \delta$ there exists a constant $C_{\delta}>0$ such that
	\[
		\| f(t, \upsilon) - f(t, \omega)\|_H \leqslant C_{\delta}\| \upsilon- \omega \|_V.
	\]
\item[(H4)] The damping coefficient $\gamma(t)\in C^{2}_{loc}([0, \infty); \mathbb{R})$ satisfies
\[
\max_{t\in \mathcal{O}}\{\gamma(t), \dot{\gamma}(t), \ddot{\gamma}(t)\}\leqslant C_{\gamma}, 
\]
where $\mathcal{O}\subset [0, \infty)$ is any compact domain and $C_{\gamma}\geqslant 0$.

\end{itemize}

By the semigroup theory (e.g., \cite[Section 6.1]{Pazy2012semigroups}) and regarding $f(t, u) - \gamma(t)v$ as a perturbation term, one can show the following lemma collecting a weak and a strong wellposedness result.
\begin{lemma}
Let Assumptions $\mathrm{(H1)}$-$\mathrm{(H4)}$ be satisfied. Then for all $u^0\in V$, $v^0\in H$, \eqref{eq:wave} has a unique weak solution
\[
	u\in C^2\left([0, T]; H^\ast \right) \cap C^1\left([0, T]; H \right) \cap C\left([0, T]; V \right)
\]
with $T = T(u^0, v^0)$. Moreover, if $u^0\in D(A)$, $v^0\in V$, there exists a strong solution to \eqref{eq:wave}
\[
	u\in C^2\left([0, T]; H \right) \cap C^1\left([0, T]; V \right) \cap C\left([0, T]; D(A) \right).
\]
\end{lemma}

We present some examples of wave-type equation that can be written in the form of \eqref{eq:secondorder}.  
First, we briefly introduce some notations. Let $\Omega $ be a bounded open region in $\mathbb{R}^{d}$, $d=2$ or $3$, with smooth boundary $\Gamma$.  We denote by $\bigtriangledown_{\Gamma}$ the gradient on $\Gamma$, and the Laplace--Beltrami operator on $\Gamma$ is given by $\triangle_{\Gamma} = \bigtriangledown_{\Gamma}\cdot \bigtriangledown_{\Gamma}$. The outer normal is denoted by $\nu$.

\begin{itemize}
\item We consider the wave equation with a homogeneous Dirichlet boundary condition
\end{itemize}
\begin{equation} \label{eq:rbc}
\begin{aligned}
            u_{tt} + \gamma(t)u_t - c_{\Omega}\triangle u &= f(t, u),  \quad \textrm{in $\mathbb{R}^{+} \times \Omega$,} \\
                u &= 0,  \quad\textrm{on $\mathbb{R}^{+} \times\Gamma$,}\\
                                        u(0,x)= u_{0}(x), \quad u_t(0, x) &= u_{1}(x),  \quad\textrm{in $\overline{\Omega}$.}
\end{aligned}
\end{equation}
For the corresponding abstract formulation of \eqref{eq:rbc}, we choose $V= H^1_0(\Omega)$ and $H = L^2(\Omega)$. The product $m(\cdot, \cdot)$ is the standard $L^2(\Omega)$ inner product, the bilinear form $a(\cdot, \cdot)$ is given by $a(\upsilon, \omega) = c_{\Omega}\int_{\Omega}\bigtriangledown\upsilon\cdot\bigtriangledown\omega\mathrm{d}x$ and $b(\cdot, \cdot)$ equals to zero. Finally, the domain of $A$ is given by $D(A) = H^1_0(\Omega)\bigcap H^2(\Omega)$.

%\item The wave equations with Robin boundary conditions.
%\end{itemize}
%\begin{equation} \label{eq:rbc}
%\begin{aligned}
%            u_{tt} + \gamma(t)u_t - c_{\Omega}\triangle u &= f_1(u),  \quad \textrm{in $\mathbb{R}^{+} \times \Omega$,} \\
%                c_{\Omega} \partial_{\nu}u + u_t &= f_2(u),  \quad\textrm{on $\mathbb{R}^{+} \times\Gamma$,}\\
%                                        u(0,x)= u_{0}(x), \quad u_t(0, x) &= u_{1}(x),  \quad\textrm{in $\overline{\Omega}$.}
%\end{aligned}
%\end{equation}

\begin{itemize}
\item The wave equation with kinetic boundary condition is considered 
\end{itemize}
\begin{equation} \label{eq:kbc}
\begin{aligned}
            u_{tt} + \left[\alpha_{\Omega} + \beta_{\Omega}\cdot\bigtriangledown + \gamma(t) \right]u_t -c_{\Omega} \triangle u &= f_1(t, u),  \quad \textrm{in $\mathbb{R}^{+} \times \Omega$,} \\
                 \mu u_{tt} + (\alpha_{\Gamma} + \beta_{\Gamma}\cdot\bigtriangledown_{\Gamma})u_t - c_{\Gamma}\triangle_{\Gamma}u &= -c_{\Omega}\partial_{\nu}u + f_2(t, u),  \quad\textrm{on $\mathbb{R}^{+} \times\Gamma$,}\\
                                        u(0,x)= u_{0}(x), \quad u_t(0, x) &= u_{1}(x),  \quad\textrm{in $\overline{\Omega}$.}
\end{aligned}
\end{equation}
To write \eqref{eq:kbc} as the abstract wave equation, we set $V= \{v\in H^1(\Omega): v|_{\Gamma}\in H^1(\Gamma)\}$ and $H = L^2(\Omega)\times L^2(\Gamma)$. The domain of $A$ is equal to the space $\{v\in H^2(\Omega): v|_{\Gamma}\in H^2(\Gamma)\}$. The inner product $m(\cdot, \cdot)$ is given by
\[
	m((v, \upsilon), (w, \omega))=\int_{\Omega}vw\mathrm{d}x + \mu \int_{\Gamma}\upsilon\omega\mathrm{d}\Gamma \quad \textrm{for any $(v, \upsilon), (w, \omega)\in H$}
\] 
and the bilinear forms
\begin{equation*}
\begin{split}
	a(\upsilon, \omega) &=c_{\Omega} \int_{\Omega}\bigtriangledown\upsilon\cdot\bigtriangledown\omega\mathrm{d}x + c_{\Gamma}\int_{\Gamma}\bigtriangledown_{\Gamma}\upsilon\cdot\bigtriangledown_{\Gamma}\omega\mathrm{d}\Gamma \quad \textrm{for any $\upsilon \in V$, $\omega\in V$,}\\
	b(\upsilon, \omega) & =\int_{\Omega}(\alpha_{\Omega}\upsilon + \beta_{\Omega}\cdot\bigtriangledown\upsilon)\omega\mathrm{d}x + \int_{\Gamma}(\alpha_{\Gamma}\upsilon + \beta_{\Gamma}\cdot\bigtriangledown_{\Gamma}\upsilon)\omega\mathrm{d}\Gamma \quad \textrm{for any $\upsilon \in V$, $\omega\in V$.}
\end{split}
\end{equation*}

\begin{itemize}
\item Consider the wave equations with acoustic boundary conditions
\end{itemize}
\begin{equation} \label{eq:abc}
\begin{aligned}
            u_{tt} + \gamma(t)u_t - c_{\Omega}\triangle u &= f_1(t, u), \quad\textrm{in $\mathbb{R}^{+} \times \Omega$,} \\
                 \mu  w_{tt} + d w_t + k_{\Gamma}w - c_{\Gamma} \triangle_{\Gamma}w &= -c_{\Omega} u_t + f_2(t, u), \quad\textrm{on $\mathbb{R}^{+} \times\Gamma$,}\\
                 \partial_{\nu}u  &= w_t, \quad\textrm{on $\mathbb{R}^{+} \times\Gamma$,}\\
                                        u(0,x)= u_{0}(x), \quad u_t(0, x) = u_{1}(x), \quad w(0,x) &= w_{0}(x), \quad w_t(0, x) = w_{1}(x), \quad\textrm{in $\overline{\Omega}$.}
\end{aligned}
\end{equation}
We use the product spaces $V= H^1(\Omega)\times H^1(\Gamma)$, $H = L^2(\Omega)\times L^2(\Gamma)$ and $D(A) \subseteq  H^2(\Omega)\times H^2(\Gamma)$. For any $\mathbf{v}=(v, \upsilon)\in H$, $\mathbf{w} = (w, \omega)\in H$, the inner product and and the bilinear forms are defined via
\begin{equation*}
\begin{split}
	m(\mathbf{v}, \mathbf{w})&=\int_{\Omega}vw\mathrm{d}x + \mu \int_{\Gamma}\upsilon\omega\mathrm{d}\Gamma,\\
	a(\mathbf{v}, \mathbf{w}) &=c_{\Omega} \int_{\Omega}\bigtriangledown v\cdot\bigtriangledown w\mathrm{d}x + \int_{\Gamma}(k_{\Gamma}\upsilon\omega+c_{\Gamma}\bigtriangledown_{\Gamma}\upsilon\cdot\bigtriangledown_{\Gamma}\omega)\mathrm{d}\Gamma, \\ 
	b(\mathbf{v}, \mathbf{w}) &= c_{\Omega}\int_{\Gamma}(v\omega-\upsilon w )\mathrm{d}\Gamma.
\end{split}
\end{equation*}
Then we can write \eqref{eq:abc} as the abstract formulation. 

%%%%%%%%%%%%%%%%%%%%%%%%%%%%%%%%%%%%%%%%%%%%%%%%%%%%%%%%%%%%%%%%%%%%%%%%%%%%%%%%

\section{Abstract full discretization}
\label{sec:full_discretization}
In this section, we integrate a vanilla or revised IMEX time scheme with an abstract space discretization to formulate a fully discrete scheme. 

\subsection{Abstract space discretization}
\label{sec:space_discretization}
Let $(V_h)_h$ be a family of finite dimensional vector spaces for the spatial approximation with respect to a discretization parameter $h$. Then we consider the semidiscretized equation of \eqref{eq:secondorder}
\begin{equation} \label{eq:spatial_discrete_secondorder}
\begin{aligned}
	\langle\ddot{u}_h, \phi_h \rangle_{V_h} + b_h(\dot{u}_h, \phi_h) + \gamma(t) m_h(\dot{u}_h, \phi_h) +a_h(u_h, \phi_h) &= \langle f_h(t, u_h), \phi_h \rangle_{V_h}  \quad t\in(0, T], \\
                                        u_h(0)= u_h^{0}, \quad \dot{u}_h(0) &= v_h^0.
\end{aligned}
\end{equation}
for all $\phi_h \in V_h$. Here, $a_h$, $b_h$ $f_h$, $m_h$, $\mathrm{p}_h$, $u_h^{0}$ and $v_h^0$ are the approximations of their corresponding continuous counterparts and share the similar properties as in Assumptions (H1)-(H3). The corresponding constants used in the discrete case are denoted by $\hat{\alpha}$, $\hat{\beta}$, $\hat{C}_{\mathrm{em}}$, $\hat{C}_{\mathrm{b}}$ and $\hat{C}_{\delta}$.

The semidiscretized equation \eqref{eq:spatial_discrete_secondorder} is equivalent to
\begin{equation}\label{eq:spatial_discrete_wave}
\begin{aligned}
\ddot{u}_{h}(t) + B_{h}\dot{u}_{h}(t) + \gamma(t)\dot{u}_{h}(t) + A_{h}u_{h}(t) &= f_{h}(t, u_{h}(t)), \quad t\geqslant 0,\\
u_{h}(0) =u_{h}^0, \quad \dot{u}_{h}(0) &= v_{h}^0.
\end{aligned}
\end{equation}
Analogously to the continuous case in Section \ref{sec:continuous_case}, we set $\mathbf{u}_h = [u_h, v_h]^\top$, $\mathbf{u}^0 = [u^0_h, v^0_h]^\top$ and  
\begin{equation*}
S_h =\left[
\begin{array}{cc}
0 & -I\\
A_h & B_h
\end{array}
\right],
\quad
g_h(t, \mathbf{u}_h)=
\left[
\begin{array}{c}
0\\
f_h(t, u_h) - \gamma(t)v_h
\end{array}
\right].
\end{equation*}
Let $H_h$ be the space $V_h$ equipped with the scalar product $m_h(\cdot, \cdot)$. Then we have the spatial discretization for the first-order evolution equation \eqref{eq:evolution} in $X_h = V_h \times H_h$ as follows.
\begin{equation}\label{eq:spatial_discrete_evolution}
\begin{aligned}
\dot{\mathbf{u}}_h + S_h\mathbf{u}_h &= g_h(t, \mathbf{u}_h), \quad t\geqslant 0,\\
\mathbf{u}_h(0) &= \mathbf{u}_h^0.
\end{aligned}
\end{equation}

In the numerical experiments (see Section \ref{sec:simulation}), we utilize the bulk-surface finite element method, as initially proposed in reference~\cite{Elliott2013}, to discretize semilinear wave equations within space $\Omega$. It is worth noting that the computational domain~$\Omega_h$ is not identical to $\Omega$ but serves as an approximation. Consequently, the finite element functions in $X_h$ are defined in $\Omega_h$, not in $\Omega$.
Then we have $X_h\nsubseteq X$\footnote{Space discretization in this situation is called nonconforming, see \cite[Definitaion 2.5]{Hipp2019}.}, which complicates direct comparison between the approximation $\mathbf{u}_h\in X_h$ and the true solution $\mathbf{u} \in X$.
To conduct an error analysis, some suitable operators are required to be a bridge to link the approximation $\mathbf{u}_h$ and the solution $\mathbf{u}$.  
\begin{itemize}
\item A lift operator $\mathcal{L}^V_h: V_h\rightarrow V$ is a bounded linear operator satisfying
	\[
		\|\mathcal{L}^V_h \upsilon_h\|_H \leqslant C_{\mathrm{H}}\| \upsilon_h\|_{H_h}, \quad \|\mathcal{L}^V_h \upsilon_h\|_V \leqslant C_{\mathrm{V}}\| \upsilon_h\|_{V_h},  \quad \forall \upsilon_h \in V_h
	\]
	with constants $C_{\mathrm{V}}, C_{\mathrm{H}} > 0$, and the adjoint lift operator $\mathcal{L}^{H^\ast}_h: H\rightarrow V_h$ and $\mathcal{L}^{V^\ast}_h: V\rightarrow V_h$ are defined via
	\begin{equation*}
	\begin{split}
		m_h(\mathcal{L}^{H^\ast}_h \upsilon, \omega_h) &= m(\upsilon, \mathcal{L}^{V}_h \omega_h), \quad \forall\upsilon \in H, \omega_h \in H_h\\
		\mathrm{p}_h(\mathcal{L}^{V^\ast}_h \upsilon, \omega_h) &= \mathrm{p}(\upsilon, \mathcal{L}^{V}_h \omega_h), \quad \forall\upsilon \in V, \omega_h \in V_h.
	\end{split}
	\end{equation*}
\item An interpolation operator $\mathcal{I}_h$ is a bounded linear operator from a dense subspace $Z^V$ of $V$ to $V_h$. 
%satisfying
%	\begin{equation*}
%		\|\mathcal{I}_h\| \leqslant C_{\mathrm{Ip}}
%	\end{equation*}
%\item We define the reference operator
\end{itemize}

\subsection{Full discretization with vanilla IMEX scheme}
\label{sec:vanilla_IMEX}

Let us denote by $\tau>0$ the step size of the time discretization and set $t_n=n\tau$, $n\geqslant 0$.
The full discrete approximations of $u$ and $\dot{u}$ are denoted by $u_h^n= u_h(t_n) \approx u(t_n)$ and $v_h^n= v_h(t_n) \approx \dot{u}(t_n)$, respectively. To simplify the following presentation, we use the short notation 
\[
	v_h^{n+\frac{1}{2}} =\frac{1}{2}(v_h^{n+1} + v_h^n), \quad f_h^n = f_h(t_n, u_h^n), \quad \gamma^n = \gamma(t_n).
\]
We apply the vanilla IMEX time discrete scheme \eqref{IMEX_1}-\eqref{IMEX_3} in~\ref{sec:IMEX} to~\eqref{eq:spatial_discrete_evolution}, and obtain the following full discrete scheme
\begin{equation} \label{full_discretization}
\begin{aligned}
v_h^{n+\frac{1}{2}} &= v_h^n - \frac{\tau}{2} A_h u_h^{n} - \frac{\tau^2}{4}A_hv_h^{n+\frac{1}{2}} - \frac{\tau}{2} B_h v_h^{n+\frac{1}{2}} - \frac{\tau}{2}\gamma^{n+1} v_h^{n+\frac{1}{2}} + \frac{\tau}{4} (\gamma^{n+1} - \gamma^{n})v_h^{n} + \frac{\tau}{2}f_h^{n},\\
u_h^{n+1} &= u_h^n + \tau v_h^{n+\frac{1}{2}},\\
v_h^{n+1} &= v_h^{n+\frac{1}{2}}  - \frac{\tau}{2} A_h u_h^{n} - \frac{\tau^2}{4}A_h v_h^{n+\frac{1}{2}} - \frac{\tau}{2} B_h v_h^{n+\frac{1}{2}} - \frac{\tau}{2}\gamma^{n+1} v_h^{n+\frac{1}{2}} + \frac{\tau}{4} (\gamma^{n+1} - \gamma^{n})v_h^{n}  + \frac{\tau}{2} f_h^{n+1}.
\end{aligned}
\end{equation}
%Here, $v_h^{n+\frac{1}{2}}: =\frac{1}{2}(v_h^{n+1} + v_h^n)$.

To derive an error bounds for this full discrete scheme, we define the data, the interpolation and the conformity error as follows
\begin{equation*}\label{space_error}
\small
\begin{split}
\epsilon^{\mathrm{data}}_h &:= \left\|u^0_h - \mathcal{L}^{V^\ast}_h u^0 \right\|_{V_h} + \left\|v^0_h - \mathcal{I}_h v^0 \right\|_{V_h} + \left\|\mathcal{L}^{H^\ast}_h f(u) - f_h(\mathcal{L}^{V^\ast}_h u) \right\|_{L^{\infty}([0, T]; H_h)},\\
\epsilon^{\mathrm{ip}}_h &:= \left\|(\mathcal{L}^{V}_h\circ\mathcal{I}_h - I)u \right\|_{L^{\infty}([0, T]; V)} + \left\|(\mathcal{L}^{V}_h\circ\mathcal{I}_h - I) \dot{u} \right\|_{L^{\infty}([0, T]; V)} + \left\|(\mathcal{L}^{V}_h\circ\mathcal{I}_h - I) \ddot{u} \right\|_{L^{\infty}([0, T]; H)},\\
\epsilon^{\mathrm{p}}_h &:= \left\|\max_{\|\varphi_h\|_{V_h} = 1}\left\{|\mathrm{p}(\mathcal{L}^{V}_h\circ\mathcal{I}_h u, \mathcal{L}^{V}_h\varphi_h) - \mathrm{p}_h(\mathcal{I}_h u,  \varphi_h)| + |\mathrm{p}(\mathcal{L}^{V}_h\circ\mathcal{I}_h \dot{u}, \mathcal{L}^{V}_h\varphi_h) - \mathrm{p}_h(\mathcal{I}_h \dot{u},  \varphi_h)|\right\} \right\|_{L^{\infty}[0, T]},\\
&\quad +\left\|\max_{\|\phi_h\|_{H_h} = 1}\left\{|m(\mathcal{L}^{V}_h\circ\mathcal{I}_h u, \mathcal{L}^{V}_h\phi_h) - m_h(\mathcal{I}_h u,  \phi_h)| + |m(\mathcal{L}^{V}_h\circ\mathcal{I}_h \ddot{u}, \mathcal{L}^{V}_h\phi_h) - m_h(\mathcal{I}_h \ddot{u},  \phi_h)|\right\} \right\|_{L^{\infty}[0, T]},\\
\epsilon^{\mathrm{b}}_h &:= \left\|\max_{\|\phi_h\|_{H_h} = 1}|b(\dot{u}, \mathcal{L}^{V}_h \phi_h) - b(\mathcal{I}_h\dot{u},  \phi_h)| \right\|_{L^{\infty}[0, T]}.
\end{split}
\end{equation*}
Based on these preparations above, we now give the statement of our main result.	
\begin{theorem}
\label{main_result}
Suppose that there exist a lift operator $\mathcal{L}^V_h$ and an interpolation operator $\mathcal{I}_h$ defined as above. 
Let $u\in C^4\left([0, T]; H \right) \cap C^3\left([0, T]; V \right) \cap C^2\left([0, T]; D(A) \right)$ be the solution of equation \eqref{eq:wave} with $u$, $\dot{u}$, $\ddot{u}\in L^{\infty}\left( [0, T]; Z^V\right)$.
If we set
\[
	\delta = \max\left\{ \max_{t_n \leqslant T}\|u_h^n\|_{V_h}, C_{\mathrm{V}}\|u\|_{L^{\infty}([0, T]; V)}\right\}
\] 
and the step size $\tau$ satisfies
\begin{equation}
\label{condition_timestep}
	\max\left\{\frac{\tau}{2}\left(\hat{ \alpha}\frac{\hat{C}_{\mathrm{em}}}{2} + \hat{\beta} \right), \frac{\tau^2}{2}\hat{\alpha} + \tau\hat{\beta}, \tau\left((1+\sqrt{3})\hat{C}_{\delta} + C_{\gamma} + \frac{3\sqrt{2}}{4}\hat{C}_{\delta}C_{\gamma} \right), \tau \right\} < 1,
\end{equation}
then for any $n$ with $t_n < T$ the full discrete approximation $u_h^{n}$, $v_h^{n}$ obtained by the scheme \eqref{full_discretization} satisfy the error bound
\[
	\|\mathcal{L}^V_h u_h^n - u(t_n)\|_V + \|\mathcal{L}^V_h v_h^n - \dot{u}(t_n)\|_H\leqslant Ce^{M t_n}\left( \epsilon^{\mathrm{total}}_h + \tau^2 + C_{\gamma}\tau \right)
\]
with $\epsilon^{\mathrm{total}}_h := \epsilon^{\mathrm{data}}_h + \epsilon^{\mathrm{ip}}_h + \epsilon^{\hat{\mathrm{a}}}_h + \epsilon^{\mathrm{b}}_h$.
\end{theorem}

%\begin{corollary}\label{main_coro}
%If $\epsilon^{\mathrm{total}}_h \leqslant h^k$ for some $k > 0$, then we obtain
%\[
%	\max_{t_n < T}\left\{\|\mathcal{L}^V_h u_h^n - u(t_n)\|_V + \|\mathcal{L}^V_h v_h^n - \dot{u}(t_n)\|_H \right\}\leqslant Ce^{M T}( h^k + \tau^2 + C_{\gamma}\tau ).
%\]
%\end{corollary}

Our main theorem will be proved in Section~\ref{sec:proof}. When the coefficient $\gamma(t) = 0$ which implies $C_\gamma = 0$, then Theorem \ref{main_result} will be reduced to \cite[Theorem 3.3]{Hochbruck2021imex}. In such a situation, the error of the scheme \eqref{full_discretization} is of oder two with respect to the time-step size. However, due to the time-dependent damping term, it is clearly from Theorem \ref{main_result} that the scheme \eqref{full_discretization} for the wave-type equations~\eqref{eq:wave} converges with order one until the error of the space discretization is reached, which is demonstrated numerically in Section~\ref{sec:simulation}.

%%%%%%%%%%%%%%%%%%%%%%%%%%%%%%%%%%%%%%%%%%%%%%%%%%%%%%%%%%%%%%%%%%%%%%%%%%%%%%%%
\subsection{Full discretization with revised IMEX scheme}
\label{sec:revised_full_discretization}
As we know, the construction of the vanilla IMEX scheme is based on two second-order methods (see~\ref{sec:IMEX}). It is reasonable to require the convergence of the IMEX scheme keep second-order for solving the semilinear wave equations with time-varying dampings. Thus, we propose a revised IMEX scheme with the abstract space discretization~\eqref{eq:spatial_discrete_evolution}, which reads 
\begin{equation} \label{full_discretization_R}
\begin{aligned}
v_h^{n+\frac{1}{2}} &= v_h^n - \frac{\tau}{2} A_h u_h^{n} - \frac{\tau^2}{4}A_hv_h^{n+\frac{1}{2}} - \frac{\tau}{2} B_h v_h^{n+\frac{1}{2}} - \frac{\tau}{2}\gamma^{n+\frac{1}{2}} v_h^{n+\frac{1}{2}} + \frac{\tau}{2}f_h^{n},\\
u_h^{n+1} &= u_h^n + \tau v_h^{n+\frac{1}{2}},\\
v_h^{n+1} &= v_h^{n+\frac{1}{2}}  - \frac{\tau}{2} A_h u_h^{n} - \frac{\tau^2}{4}A_h v_h^{n+\frac{1}{2}} - \frac{\tau}{2} B_h v_h^{n+\frac{1}{2}} - \frac{\tau}{2}\gamma^{n + \frac{1}{2}} v_h^{n+\frac{1}{2}} + \frac{\tau}{2} f_h^{n+1}.
\end{aligned}
\end{equation}
Here, the approximation of $u$ and $\dot{u}$ are still denoted by $u_h^{n}$ and $v_h^{n}$, respectively. In comparison with discretization~\eqref{full_discretization}, the term $\gamma^{n+1}$ is replaced by $\gamma^{n+\frac{1}{2}}$, and $\frac{\tau}{4}(\gamma^{n+1} -\gamma^{n})v_h^n$ is removed simultaneously.
This revised IMEX time scheme is also derived from the Crank-Nicolson scheme. Indeed, we have
%\begin{equation*}
%\begin{aligned}
% - \frac{\tau}{2}\left(\gamma^n v^{n} +\gamma^{n+1} v^{n+1}\right) &=- \frac{\tau}{2}\left(2\gamma^n v^{n+\frac{1}{2}} - \gamma^n v^{n+1}+ \gamma^{n+1}v^{n+1}\right)\\
% &= - \tau\gamma^{n} v^{n+\frac{1}{2}}  - \frac{\tau}{2} (\gamma^{n+1} - \gamma^{n})v^{n+1}
% \end{aligned}
%\end{equation*}
%\begin{equation*}
%\begin{aligned}
% - \frac{\tau}{2}\left(\gamma^n v^{n} +\gamma^{n+1} v^{n+1}\right) &=- \frac{\tau}{2}\left(\gamma^n v^{n}+2\gamma^{n+1} v^{n+\frac{1}{2}} -  \gamma^{n+1}v^{n}\right)\\
% &= - \tau\gamma^{n+1} v^{n+\frac{1}{2}}  + \frac{\tau}{2} (\gamma^{n+1} - \gamma^{n})v^{n}
% \end{aligned}
%\end{equation*}
\begin{align}
u^{n+1} &= u^n + \tau v^{n+\frac{1}{2}},\label{IMEX_R_1}\\
v^{n+1} &= v^n - \frac{\tau}{2}A\left(u^{n} + u^{n+1}\right) - \tau Bv^{n+\frac{1}{2}} - \tau\gamma^{n+1} v^{n+\frac{1}{2}} + \frac{\tau}{2} (\gamma^{n+1} - \gamma^{n})v^{n} + \frac{\tau}{2}\left(f^{n} + f^{n+1}\right),\label{IMEX_R_2}\\
v^{n+1} &= v^n - \frac{\tau}{2}A\left(u^{n} + u^{n+1}\right) - \tau Bv^{n+\frac{1}{2}} - \tau\gamma^{n} v^{n+\frac{1}{2}} - \frac{\tau}{2} (\gamma^{n+1} - \gamma^{n})v^{n+1} + \frac{\tau}{2}\left(f^{n} + f^{n+1}\right),\label{IMEX_R_3}
\end{align}
which is induced by \eqref{CN} and $2v^{n+\frac{1}{2}}  = v^{n+1} +  v^n$. By the sum of \eqref{IMEX_R_2} and \eqref{IMEX_R_3}, we get
\begin{equation} \label{IMEX_R_4}
\begin{split}
v^{n+1} = v^n - \frac{\tau}{2}A\left(u^{n} + u^{n+1}\right) - \tau Bv^{n+\frac{1}{2}} - \tau\gamma^{n+\frac{1}{2}} v^{n+\frac{1}{2}} - \frac{\tau}{4} (\gamma^{n+1} - \gamma^{n})(v^{n+1}-v^{n}) + \frac{\tau}{2}\left(f^{n} + f^{n+1}\right).
\end{split}
\end{equation}
Inserting \eqref{IMEX_R_1} into \eqref{IMEX_R_4} and eliminating $u^{n+1}$ yields
\begin{equation*}
v^{n+1} = v^n - \tau Au^{n} - \frac{\tau^2}{2}Av^{n+\frac{1}{2}} - \tau Bv^{n+\frac{1}{2}} - \tau\gamma^{n+\frac{1}{2}} v^{n+\frac{1}{2}} - \frac{\tau}{4} (\gamma^{n+1} - \gamma^{n})(v^{n+1}-v^{n})+ \frac{\tau}{2}\left(f^{n} + f^{n+1}\right).
\end{equation*}
Furthermore, with $2v^{n+\frac{1}{2}}  = v^{n+1} +  v^n$ we obtain
\begin{equation} 
\label{IMEX_R_5}
\small
\begin{split}
v^{n+\frac{1}{2}} &= v^n - \frac{\tau}{2} Au^{n} - \frac{\tau^2}{4}Av^{n+\frac{1}{2}} - \frac{\tau}{2} Bv^{n+\frac{1}{2}} - \frac{\tau}{2}\gamma^{n+\frac{1}{2}} v^{n+\frac{1}{2}} - \frac{\tau}{8} (\gamma^{n+1} - \gamma^{n})(v^{n+1}-v^{n})+ \frac{\tau}{4}\left(f^{n} + f^{n+1}\right),\\
v^{n+1} &= v^{n+\frac{1}{2}} - \frac{\tau}{2} Au^{n} - \frac{\tau^2}{4}Av^{n+\frac{1}{2}} - \frac{\tau}{2} Bv^{n+\frac{1}{2}} - \frac{\tau}{2}\gamma^{n+\frac{1}{2}} v^{n+\frac{1}{2}} - \frac{\tau}{8} (\gamma^{n+1} - \gamma^{n})(v^{n+1}-v^{n})+ \frac{\tau}{4}\left(f^{n} + f^{n+1}\right).
\end{split}
\end{equation}
We note that $ \frac{\tau}{8} (\gamma^{n+1} - \gamma^{n})(v^{n+1}-v^{n})= o(\tau^3)$. By eliminating this term and replacing the trapezoidal rule for the nonlinear part by the left/right rectangular rule, respectively in \eqref{IMEX_R_5}, we have
\begin{equation} 
\label{IMEX_R_6}
\begin{split}
v^{n+\frac{1}{2}} &= v^n - \frac{\tau}{2} Au^{n} - \frac{\tau^2}{4}Av^{n+\frac{1}{2}} - \frac{\tau}{2} Bv^{n+\frac{1}{2}} - \frac{\tau}{2}\gamma^{n+\frac{1}{2}} v^{n+\frac{1}{2}} + \frac{\tau}{2}f^{n},\\
v^{n+1} &= v^{n+\frac{1}{2}} - \frac{\tau}{2} Au^{n} - \frac{\tau^2}{4}Av^{n+\frac{1}{2}} - \frac{\tau}{2} Bv^{n+\frac{1}{2}} - \frac{\tau}{2}\gamma^{n+\frac{1}{2}} v^{n+\frac{1}{2}} + \frac{\tau}{2}f^{n+1}.
\end{split}
\end{equation}
Therefore, combing \eqref{IMEX_R_1}, \eqref{IMEX_R_6} with the abstract space discretization gives \eqref{full_discretization_R}.  
%\begin{equation*}
%\small
%\begin{split}
%v^{n+1} &= v^n - \tau Au^{n} - \frac{\tau^2}{2}A{\color{violet}v^{n+\frac{1}{2}}} - \tau Bv^{n+\frac{1}{2}} - \tau\gamma^{n+\frac{1}{2}} v^{n+\frac{1}{2}} - \frac{\tau}{4} (\gamma^{n+1} - \gamma^{n})(v^{n+1} - v^{n}) + \frac{\tau}{2}\left(f^{n} + f^{n+1}\right)\\
%	&= v^n - \tau Au^{n} - \frac{\tau}{2}A\left( u^{n+1} - u^n\right) - \tau Bv^{n+\frac{1}{2}}- \tau\gamma^{n+\frac{1}{2}} v^{n+\frac{1}{2}} - \frac{\tau}{4} (\gamma^{n+1} - \gamma^{n})(v^{n+1} - v^{n}) + \frac{\tau}{2}\left(f^{n} + f^{n+1}\right)\\
%&= v^n - \frac{\tau}{2}A\left( u^{n+1} + u^n\right) - \tau B{\color{brown}v^{n+\frac{1}{2}}} - \tau\gamma^{n+\frac{1}{2}} {\color{brown}v^{n+\frac{1}{2}}}  - \frac{\tau}{4} (\gamma^{n+1} - \gamma^{n})(v^{n+1} - v^{n}) + \frac{\tau}{2}\left(f^{n} + f^{n+1}\right)\\
%&= v^n - \frac{\tau}{2}A\left( u^{n+1} + u^n\right) - \frac{\tau}{2}B\left(v^{n} + v^{n+1}\right) - \frac{\tau}{2}\gamma^{n+\frac{1}{2}}\left(v^{n} + v^{n+1}\right) - \frac{\tau}{4} (\gamma^{n+1} - \gamma^{n})(v^{n+1} - v^{n})+ \frac{\tau}{2}\left(f^{n} + f^{n+1}\right) \\
%	&\quad - \frac{\tau^2}{4}B\left(f^{n} - f^{n+1}\right) - \frac{\tau^2}{4}\gamma^{n+\frac{1}{2}}\left(f^{n} - f^{n+1}\right)\\
%\end{split}
%\end{equation*}

In the subsequent section, we utilize this revised scheme \eqref{full_discretization_R} to solve a class of semilinear wave equations with nonautonomous dampings in the following section. It is evident that this numerical scheme achieves a convergence rate of second order in terms of time-step size. Unfortunately, we have not shown this scheme is a second-order method theoretically so far, which will be part of future research work.  

%%%%%%%%%%%%%%%%%%%%%%%%%%%%%%%%%%%%%%%%%%%%%%%%%%%%%%%%%%%%%%%%%%%%%%%%%%%%%%%%
\section{Numerical simulation}
\label{sec:simulation}
%Some numerical experiments are implemented to illustrate our main results in this Section. 
In this section, we apply the full discretization \eqref{full_discretization} and \eqref{full_discretization_R} to solve the wave equation \eqref{eq:kbc} with the damping coefficient of the form
\[
\gamma(t) = r_1(r_2 + t)^{\eta}
\]
on the unit disc $\Omega = \{x \in \mathbb{R}^2: |x|<1\}$ and its boundary $\Gamma = \{x \in \mathbb{R}^2: |x|=1\}$.

Since the boundary $\Gamma$ is smooth, it is not feasible to represent $\Omega$ precisely with a polygonal mesh. Consequently, in all subsequent numerical experiments, we employ bulk-surface finite element method, as elaborated in~\cite[Chapter 7]{Hipp2017} and~\cite{Elliott2013}, to discretize~\eqref{eq:kbc} in space.   
The smooth domain $\Omega$ is approximated by a computational bulk $\Omega_h$, which is a polygonal approximation with a triangulation. This triangulation is the mesh of isoparametric elements of degree $2$. Here, the discretization parameter $h$ represents the maximal mesh width. We refer to $\Gamma_h=\partial\Omega_h$ as the corresponding computational surface.  
The coefficients set is given in Table \ref{coefficients_set}. The code developed in this work will be made available on request.
\begin{table}[!htbp]
\begin{center}
\begin{tabular}{ccccccc}
\toprule
$\alpha_{\Omega}$ & $\beta_{\Omega}$ & $c_{\Omega}$ & $\mu$ &$\alpha_{\Gamma}$ & $\beta_{\Gamma}$ & $c_{\Gamma}$ \\
\midrule
 $0$&$0$  &$1$  &$1$ &$0$ &$0$ &$1$\\
\bottomrule
\end{tabular}
\caption{The value of the coefficients in equation \eqref{eq:kbc}.}\label{coefficients_set}
\end{center}
\end{table}

\textbf{Example 1}. We choose the nonlinearities as
\begin{equation*}
\begin{split}
	f_1(t, u) &= |u|u -\left(4\pi^2 + |\sin(2\pi t)x_1 x_2|\right)\sin(2\pi t)x_1 x_2 + 2\pi\gamma(t)\cos(2\pi t)x_1 x_2, \\
	f_2(t, u) &= u^3 -4\pi^2 \sin(2\pi t)x_1 x_2 + 2\pi\gamma(t)\sin(2\pi t)x_1 x_2 - \left( \sin(2\pi t)x_1 x_2\right)^3
\end{split}
\end{equation*}
and the initial data $u(0, x) = 0$, $\dot{u}(0, x) = 2\pi x_1 x_2$. Then we can deduce the exact solution of the equation \eqref{eq:kbc}
\[
	u(t, x) = \sin(2\pi t)x_1 x_2.
\]

To illustrate our main results, the exact solution and numerical solution are compared. Specifically, we consider the error
\[
	\mathcal{E}(t) = \|u(t, x)|_{\Omega_h} - u_h(t, x) \|_{V_h} + \| \dot{u}(t, x)|_{\Omega_h} - \dot{u}_h(t, x)\|_{H_h}
\]
with $V_h= \{v_h\in H^1(\Omega_h): v_h|_{\Gamma_h}\in H^1(\Gamma_h)\}$ and $H_h = L^2(\Omega_h)\times L^2(\Gamma_h)$.

\begin{figure}[!htbp]
\centering
\includegraphics[width=0.4\textwidth]{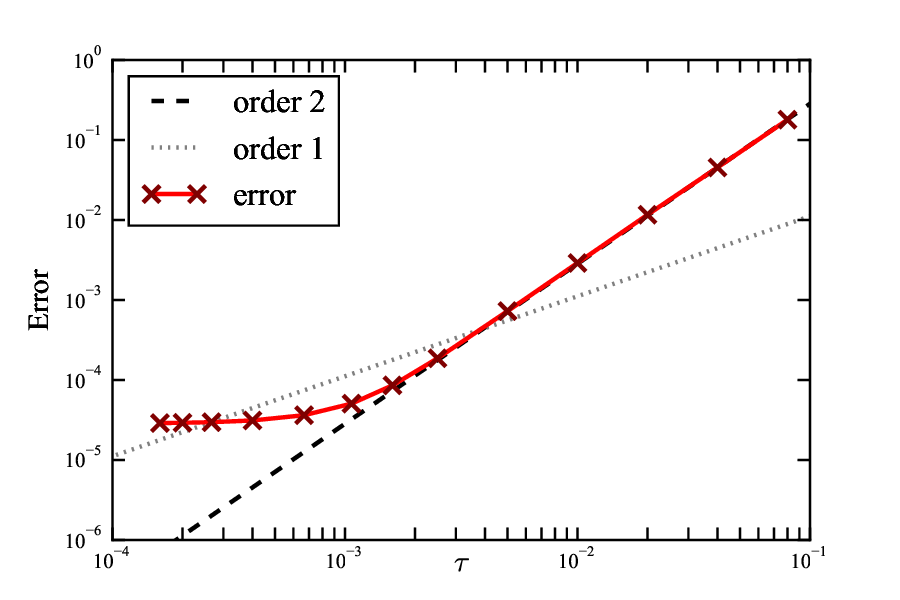}
\hspace{0.2in}
\includegraphics[width=0.4\textwidth]{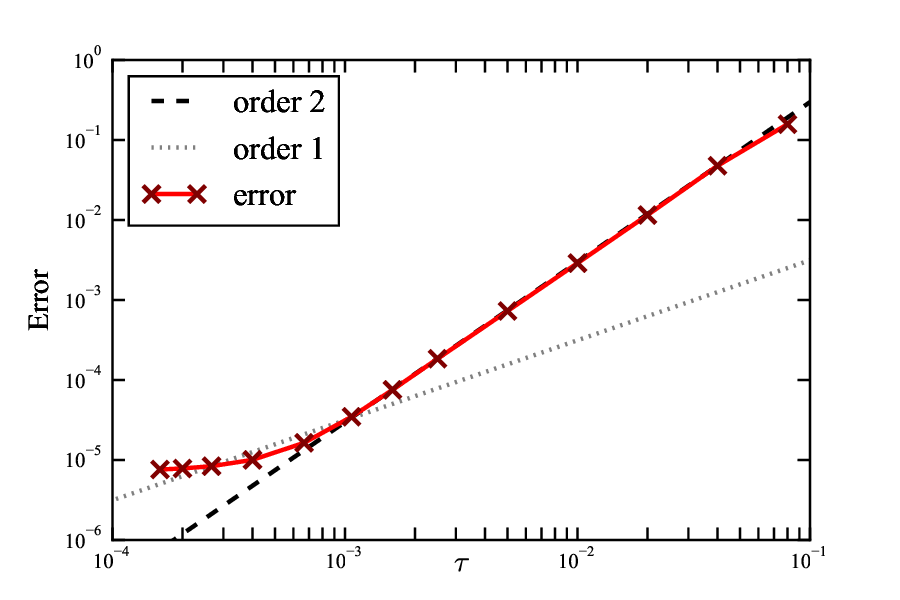}
\caption{Error $\mathcal{E}(0.8)$ of the full discretization \eqref{full_discretization} plots for the problem with the damping coefficient $\gamma(t) = 0$. Left: coarse space discretization~($328193$ degrees of freedom, $h = 0.0139201$); Right: fine space discretization ($1311745$ degrees of freedom, $h = 0.00697294$).}
\label{fig_1}
\end{figure}

If the damping coefficient $\gamma(t) = 0$ that deduces $C_{\gamma} = 0$, we can observe from Figure \ref{fig_1} (in both plots) that the temporal convergence rate of the full discretization matches $o(\tau^2)$ for a coarse and a fine space discretization, respectively. These numerical experiments are in accordance with the theoretical results in Theorem~\ref{main_result} in this paper and~\cite[Theorem 3.3]{Hochbruck2021imex}.

\begin{figure}[!htbp]
\centering
\subfigure[vanilla IMEX scheme]{
%\includegraphics[width=0.4\textwidth]{fig/fig2_coarse.eps}
%\hspace{0.2in}
\includegraphics[width=0.4\textwidth]{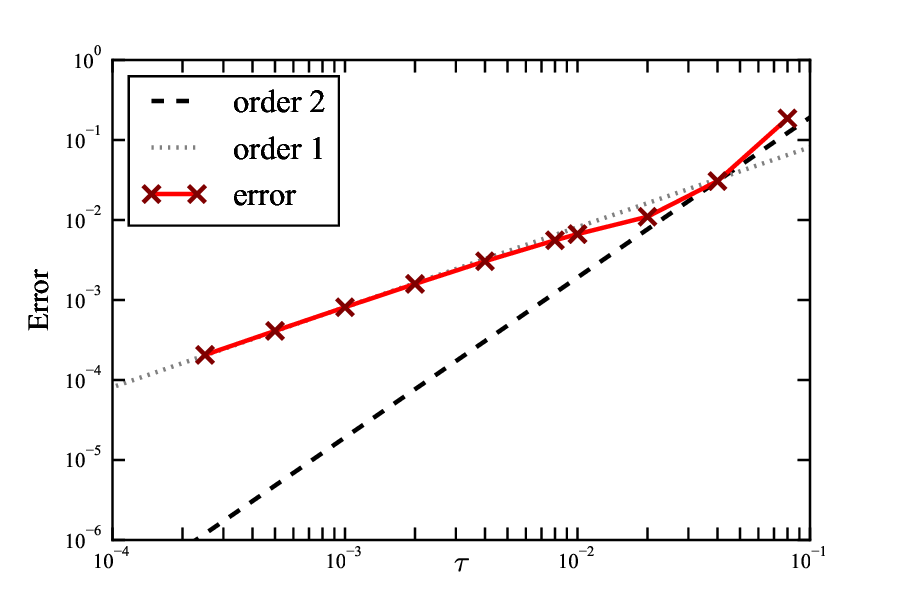}
}
\hspace{0.2in}
\subfigure[revised IMEX scheme]{
%\includegraphics[width=0.4\textwidth]{fig/fig2_coarse_R.eps}
%\hspace{0.2in}
\includegraphics[width=0.4\textwidth]{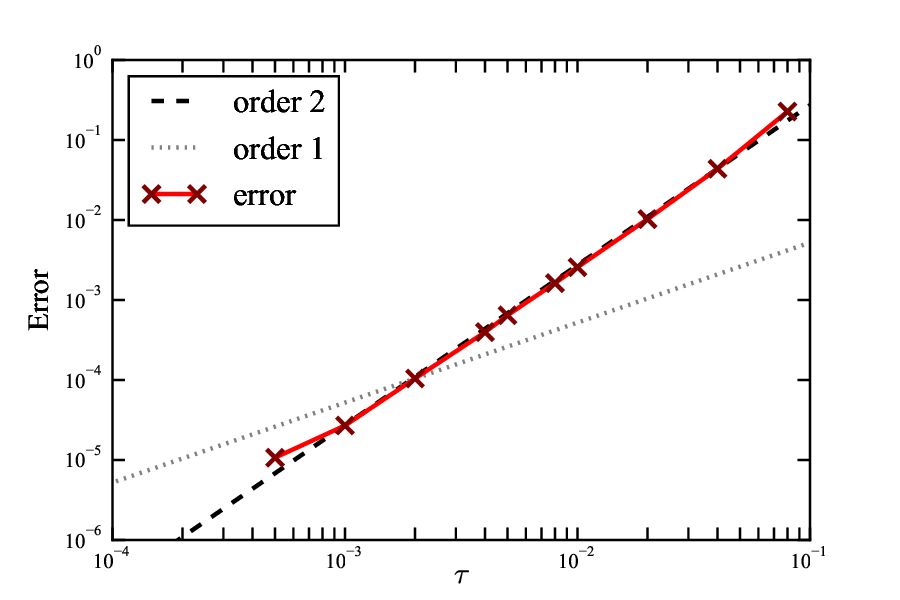}
}
\caption{Error $\mathcal{E}(0.8)$ of the full discretizations plots for the problem with the damping coefficient satisfying $r_1 =r_2 =1$ and $\eta = -2$.}
\label{fig_2}
\end{figure}

Regarding the damping coefficient $\gamma(t) = (1+t)^{-2}$, Figure~\ref{fig_2} illustrates the temporal convergence rates of the full discretizations \eqref{full_discretization} and \eqref{full_discretization_R} for the fine space discretization ($h=0.00697294$). 
It can be observed that the scheme \eqref{full_discretization} converges with order one in terms of time step in Figure~\ref{fig_2} (a),  while the scheme \eqref{full_discretization_R} converges with order two in Figure~\ref{fig_2} (b). 
The temporal convergence behaviors of both full discretizations for the problems with damping coefficients $\gamma(t) = 1+t$, and $\gamma(t) = (1+t)^2$
are very similar to those in Figure~\ref{fig_2}, which can be seen clearly in Figure~\ref{fig_4} and \ref{fig_5} of~\ref{sec:different_gamma}. 

\begin{figure}[!htbp]
\centering
\subfigure[Crank--Nicolson scheme]{
\includegraphics[width=0.4\textwidth]{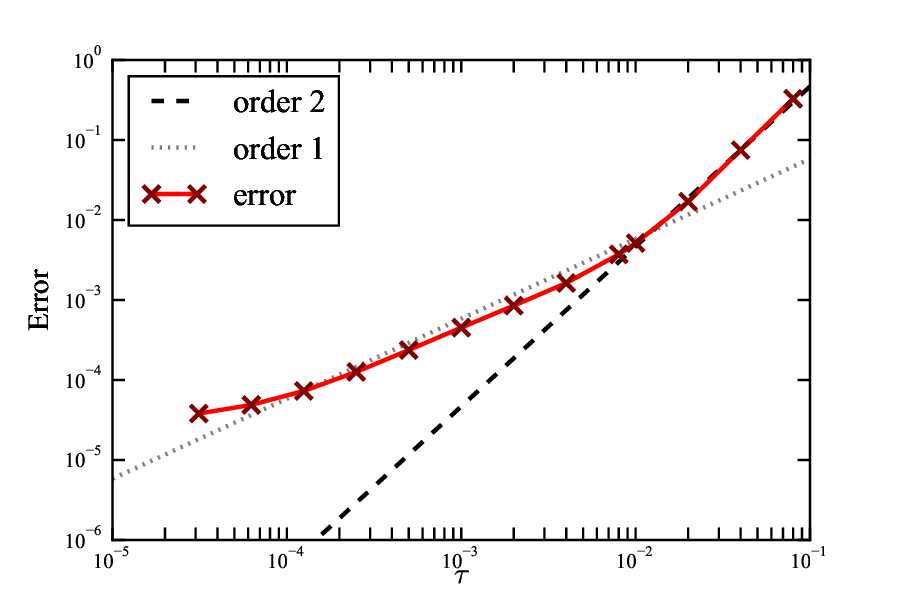}
}
\hspace{0.2in}
\subfigure[Runge--Kutta scheme]{
\includegraphics[width=0.4\textwidth]{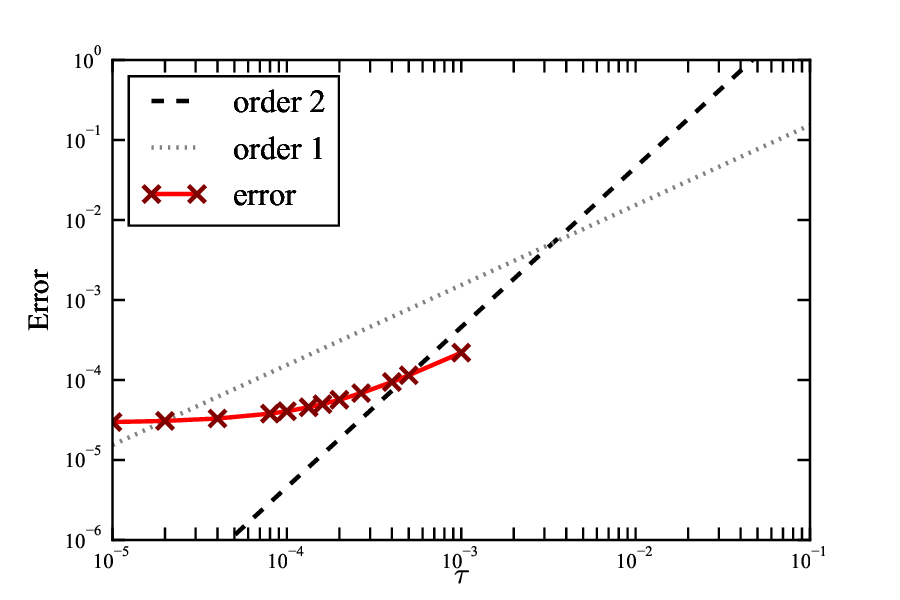}
}
\subfigure[vanilla IMEX scheme]{
%\includegraphics[width=0.4\textwidth]{fig/fig3_coarse.eps}
%\hspace{0.2in}
\includegraphics[width=0.4\textwidth]{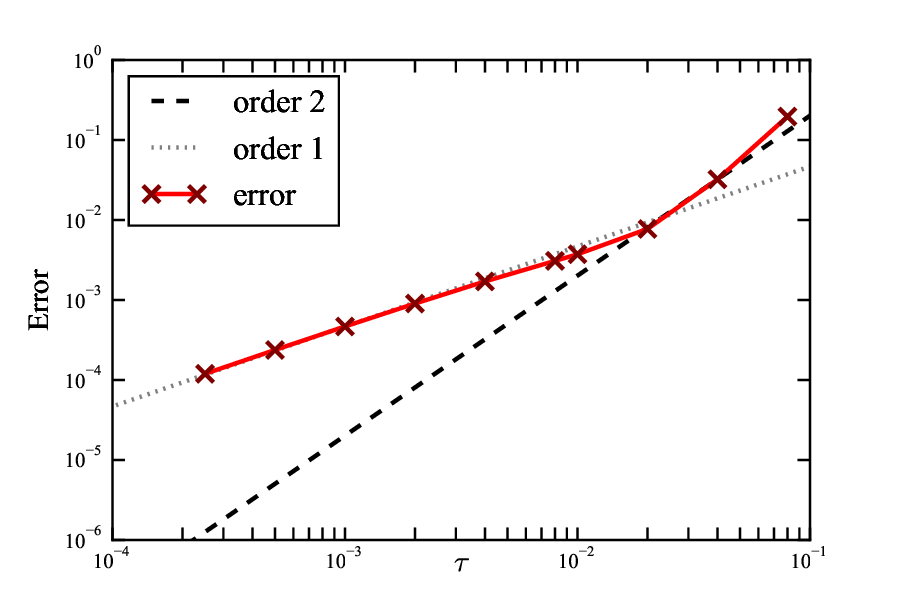}
}
\hspace{0.2in}
\subfigure[revised IMEX scheme]{
%\includegraphics[width=0.4\textwidth]{fig/fig3_coarse_R.eps}
%\hspace{0.2in}
\includegraphics[width=0.4\textwidth]{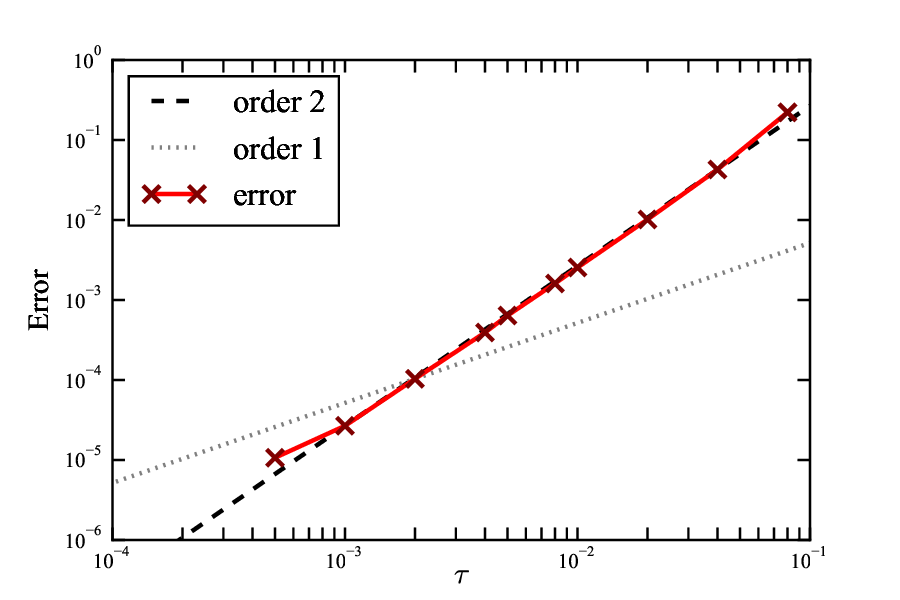}
}
\caption{Error $\mathcal{E}(0.8)$ of the full discretization plots for the problem with the damping coefficient satisfying $r_1 =r_2 = 1$ and $\eta = -1$. }
\label{fig_3}
\end{figure} 

\begin{figure}[!htbp]
\centering
\subfigure[]{
\includegraphics[width=0.6\textwidth]{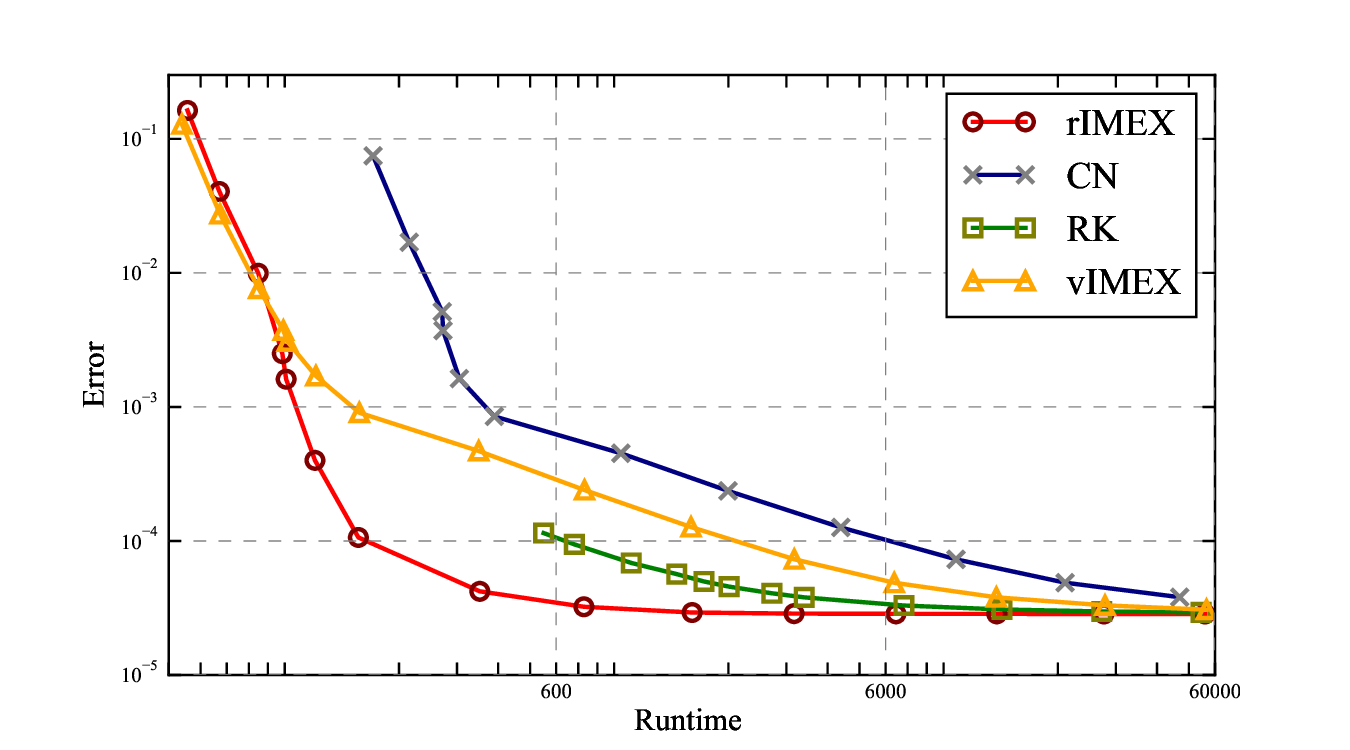}
}
\caption{Error $\mathcal{E}(0.8)$ of different time schemes with the same space discretization~($h = 0.0139201$) plotted against the runtime. Here, vIMEX: the vanilla IMEX scheme, rIMEX: the revised IMEX scheme, CN: the Crank--Nicolson scheme, and RK: the Runge--Kutta scheme.}
\label{fig_6}
\end{figure}

Finally, Figure~\ref{fig_3} shows the errors of different time discrete schemes with the coarse spatial discretization ($h = 0.0139201$) for the problem with the coefficient $\gamma(t) = (1+t)^{-1}$, while Figure~\ref{fig_6} demonstrates a comparison of these errors pictured along the runtime. As can be seen easily, the Crank--Nicolson scheme~\eqref{CN_equiv} converges with order one (see Figure~\ref{fig_3} (a)), but its convergence speed is much slower than that of  of the revised IMEX scheme~\eqref{IMEX_R_6} as in Figure~\ref{fig_6}. 
It also indicates that the error of the scheme \eqref{full_discretization_R} reaches the space discretization error plateau fast compared to the scheme \eqref{full_discretization}.
Although the classical Runge--Kutta scheme~\cite{hansen2006runge} is more efficient than the revised IMEX scheme, the Runge--Kutta scheme is strictly restricted by time step (see Figure~\ref{fig_3} (b)), which is only stable under a strong Courant--Friedrichs--Lewy condition. Additionally, Table~\ref{runtime_comparison} shows the execution time for various methods to achieve the same level of error, thereby highlighting the enhanced computational efficiency of our revised IMEX scheme.  

\begin{table}[!htbp]
\begin{center}
\begin{tabular}{cc}
\toprule
Scheme & Runtime (\textrm{s}) \\
\midrule
revised IMEX &   $351.87$ \\
vanilla IMEX &  $6375.01$\\
RK & $2009.62$ \\
CN & $21037.60$ \\
\bottomrule
\end{tabular}
\caption{Runtime of the different schemes to reach the same error $\approx 4.6\times 10^{-5}$.}\label{runtime_comparison}
\end{center}
\end{table}

\begin{figure}[!htbp]
\centering
\subfigure[$t=0$]{
\includegraphics[width=0.45\textwidth]{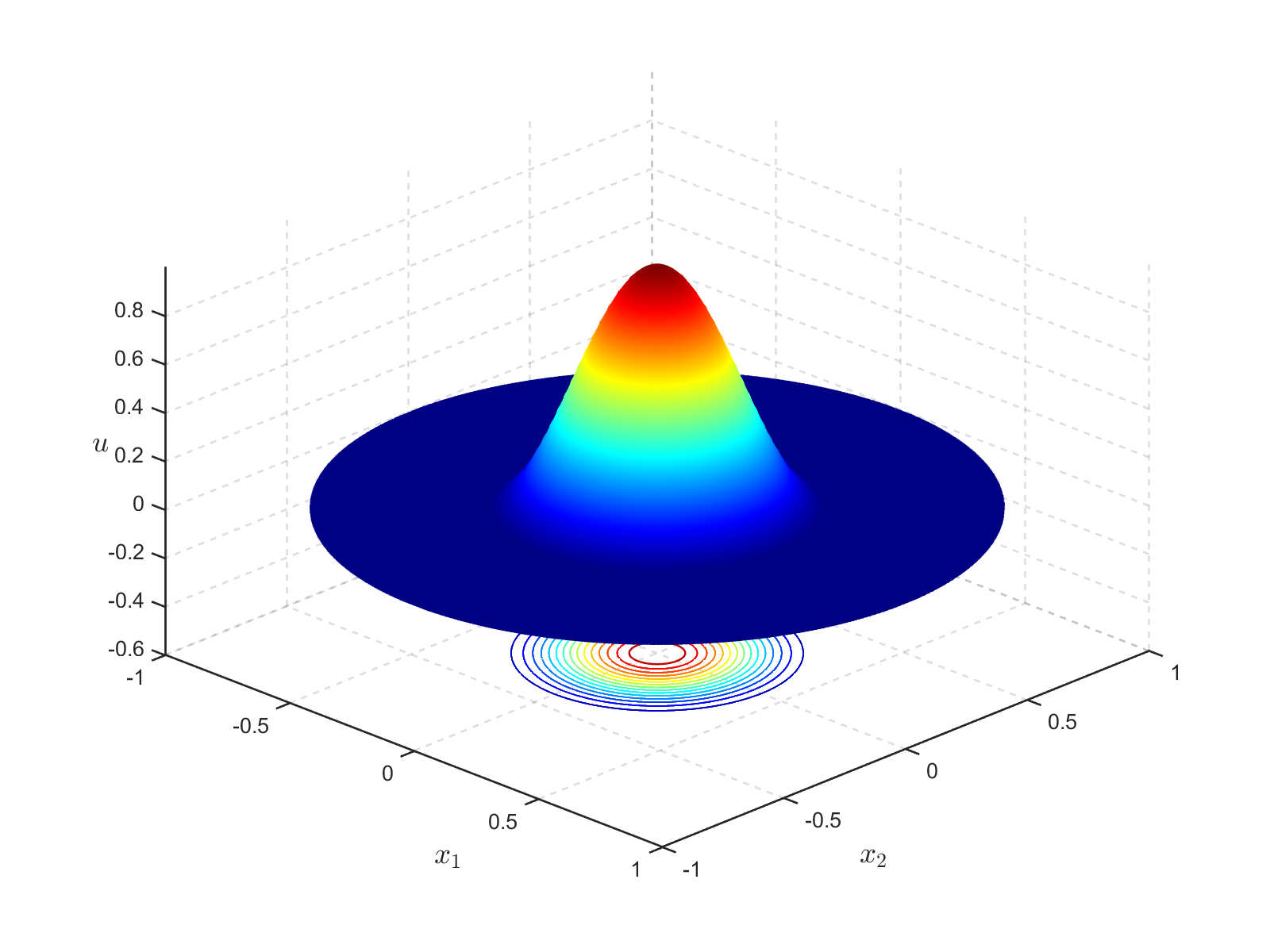}
}
%\hspace{0.2in}
\subfigure[$t=0.4$]{
\includegraphics[width=0.45\textwidth]{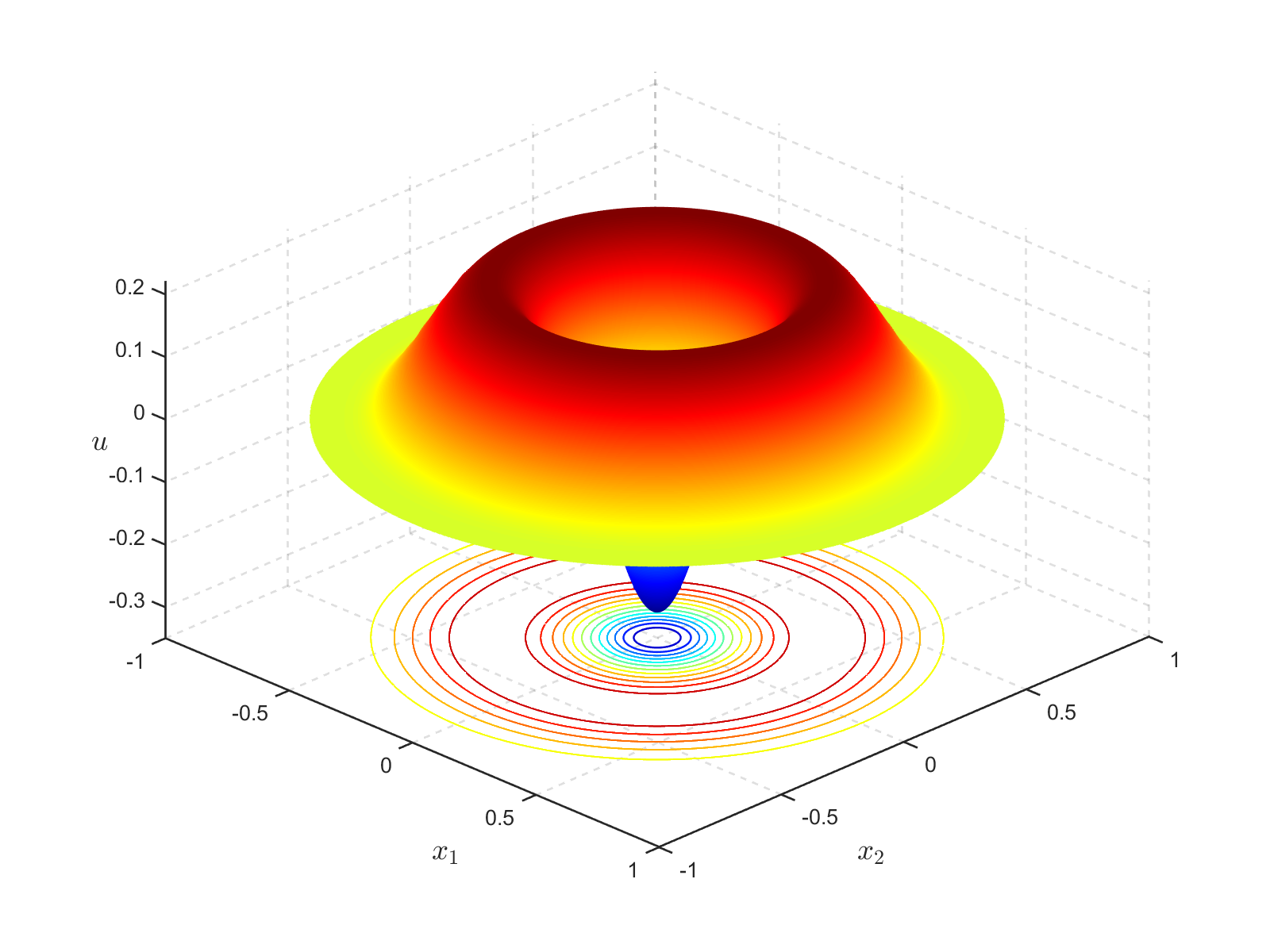}
}
\subfigure[$t=0.8$]{
\includegraphics[width=0.45\textwidth]{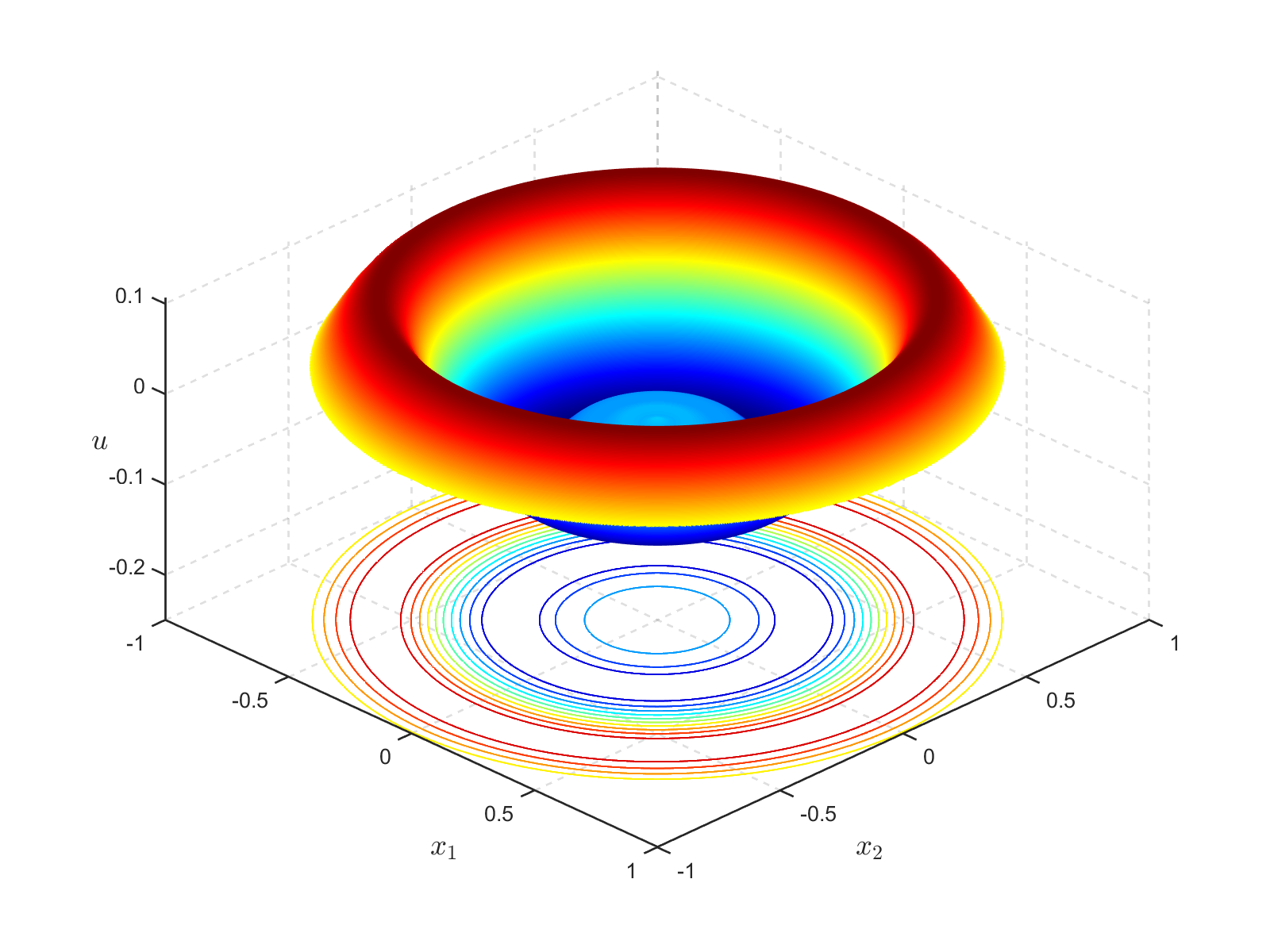}
}
%\hspace{0.2in}
\subfigure[$t=1.2$]{
\includegraphics[width=0.45\textwidth]{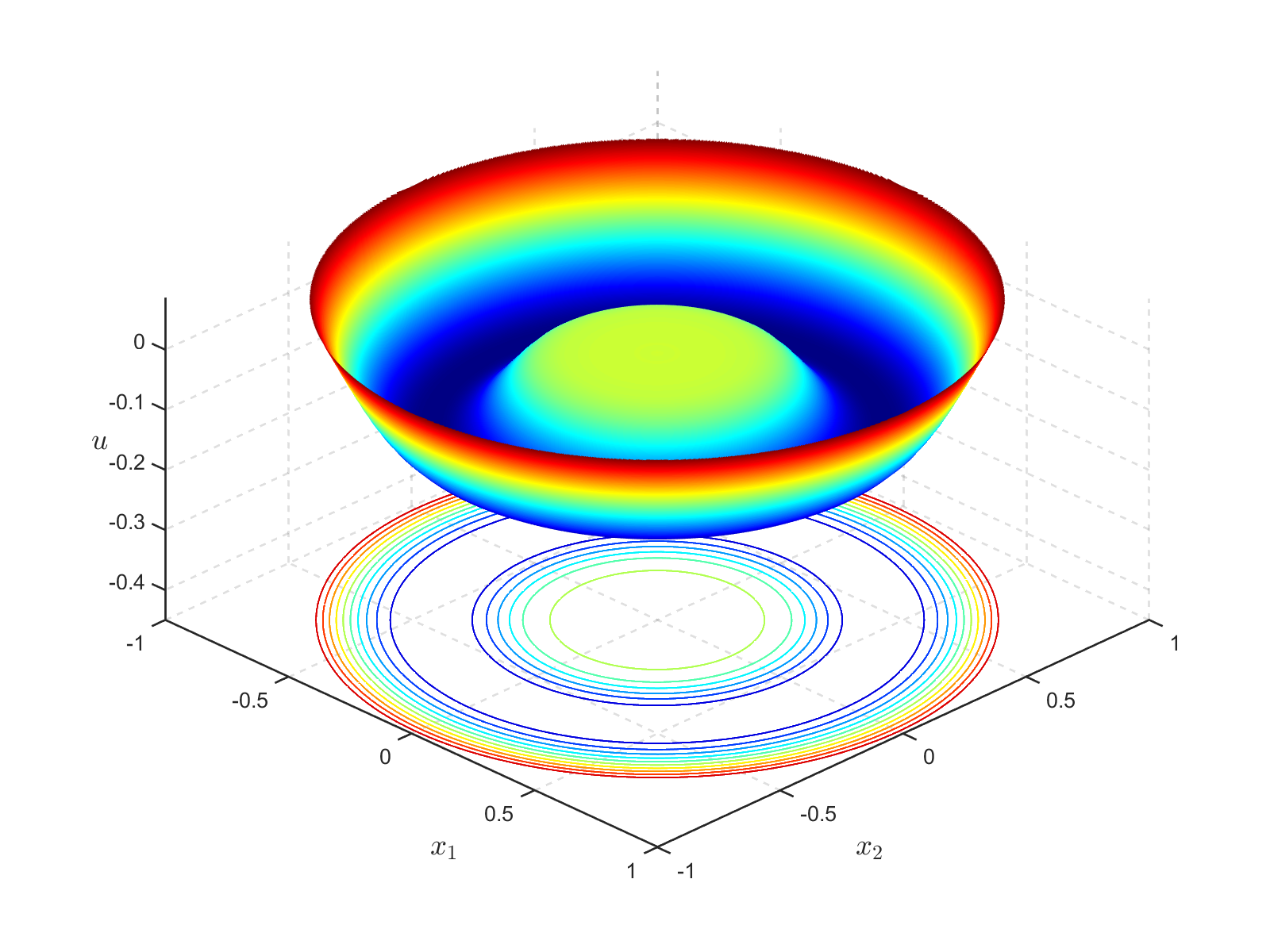}
}
%\hspace{0.2in}
\subfigure[$t=1.6$]{
\includegraphics[width=0.45\textwidth]{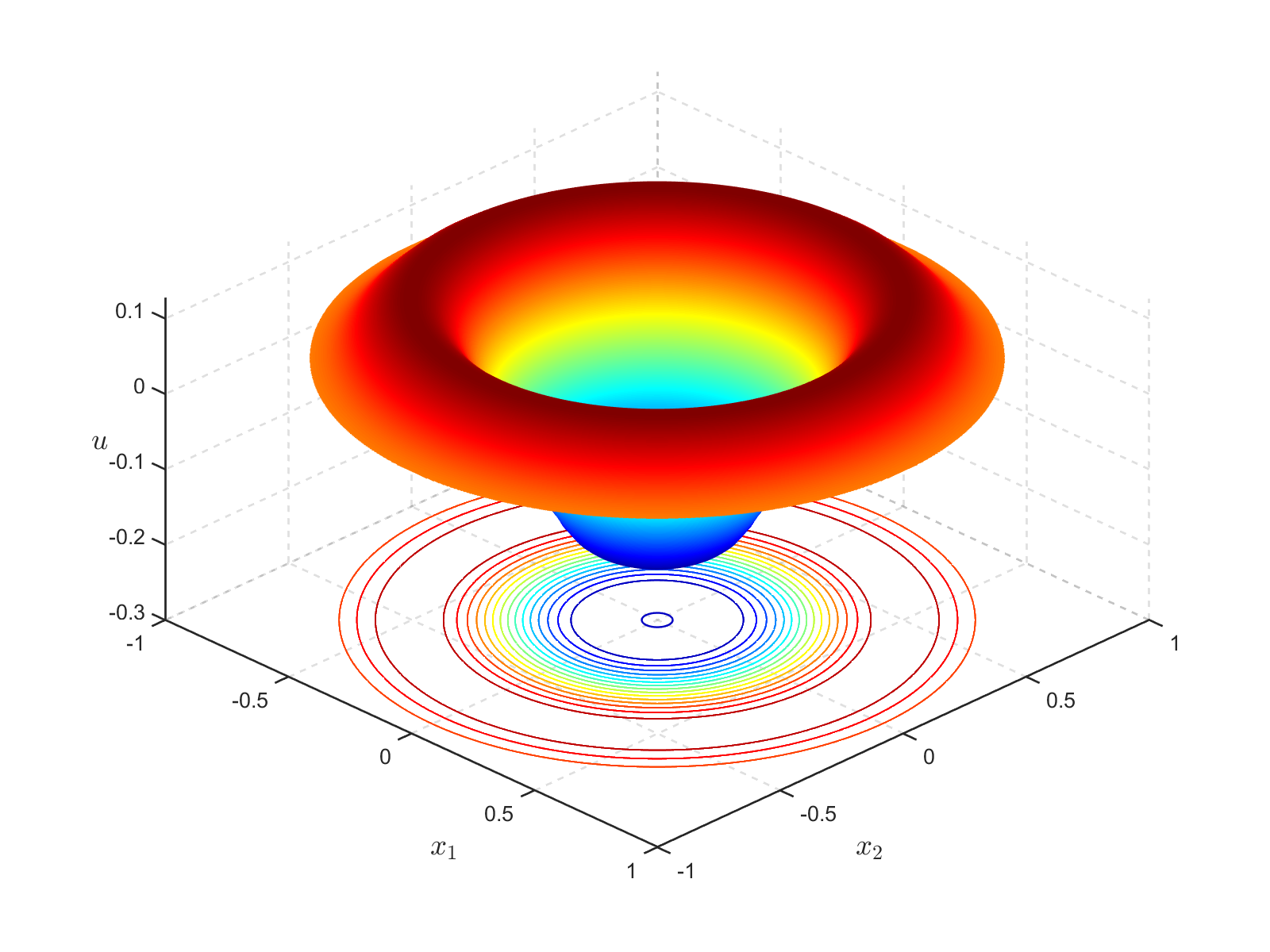}
}
%\hspace{0.2in}
\subfigure[$t=2.0$]{
\includegraphics[width=0.45\textwidth]{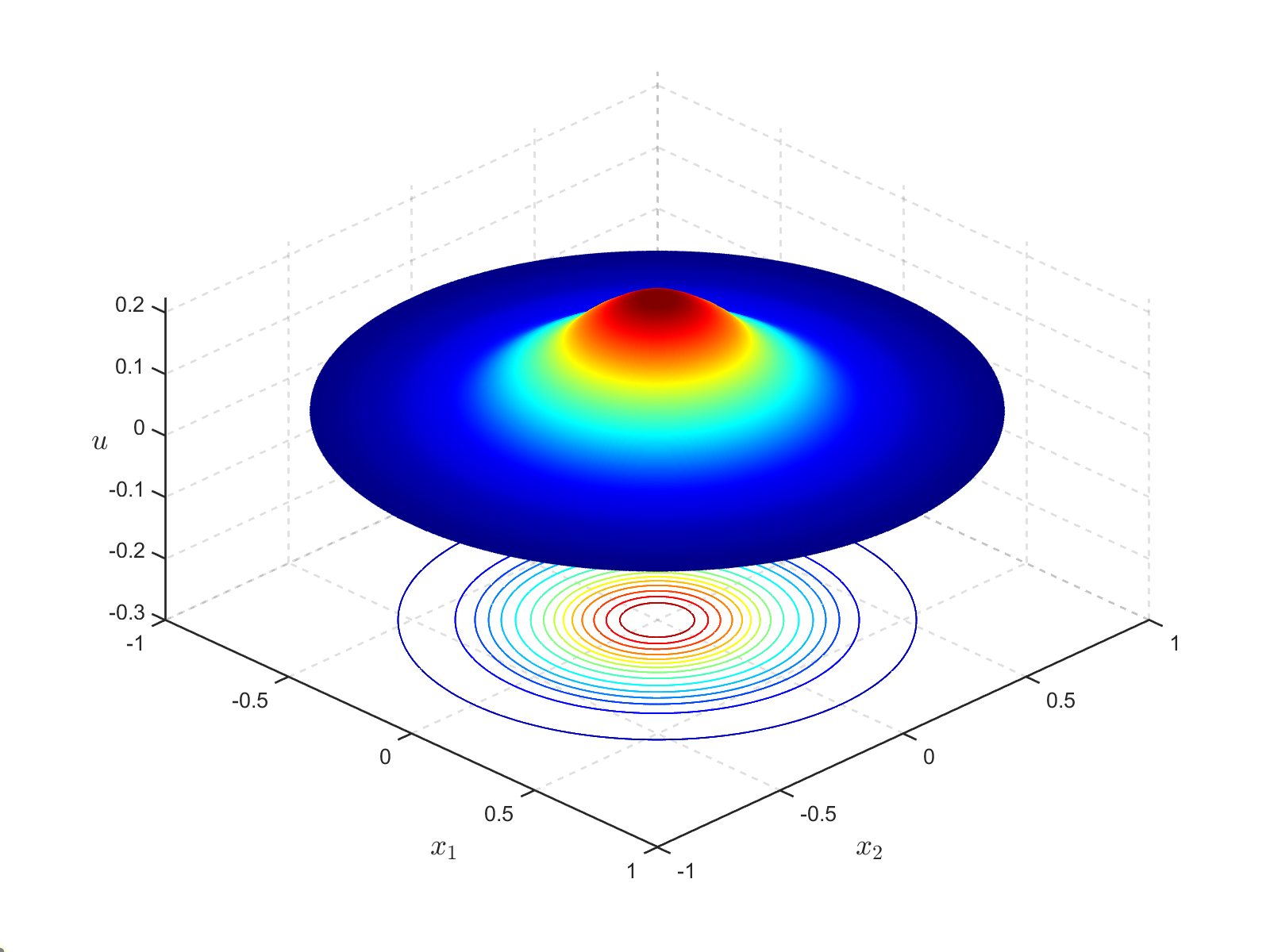}
}
\caption{The profiles of numerical solutions to the problem with the damping coefficient satisfying $r_1 =r_2 = 1$ and $\eta = -1$. }
\label{fig_7}
\end{figure}

\textbf{Example 2}. The nonlinearities is chosen as 
\begin{equation*}
	f_1(t, u) = |u|u, \quad f_2(t, u) = u^3
\end{equation*}
and the initial data is given by 
%\[
%	u(0, x) = \cos^2(2\pi\sqrt{x_1^2 + x_2^2}), \quad \dot{u}(0, x) = 0.
%\]
\begin{equation*}
u(0, x)=\left\{
\begin{array}{cc}
\cos^2\left(\pi |x|\right),& |x| \leqslant 0.5\\
0, & \textrm{otherwise}
\end{array}
\right.,
\quad
\dot{u}(0, x) = 0.
\end{equation*} 
%As we know, this problem does not have an exact solution.
Figure~\ref{fig_7} shows the propagation of the numerical solution by utilizing our revised IMEX scheme at various time points.

%%%%%%%%%%%%%%%%%%%%%%%%%%%
%%%%%%%%%%%%%%%%%%%%%%%%%%
\section{Proof of Theorem \ref{main_result}}
\label{sec:proof}
Using the similar idea (cf., e.g., \cite{Hochbruck2016error, Hochbruck2021imex}) for error analysis, we give the proof of our main theorem. Throughout this section, $C_1, C_2, \cdots$ denote corresponding constants 
\begin{proof}
We use the notation $\mathbf{u}^n_h= [u^n_h, v^n_h]^\top$ and $\tilde{\mathbf{u}}^n= [\tilde{u}^n, \tilde{v}^n]^\top = [u(t_n), \dot{u}(t_n)]^\top$ in the following. The subsequent proof consists of several steps.

\textbf{Step 1}. Note that the error of the scheme \eqref{full_discretization} can be split into two parts
\[
	\mathcal{L}_h\mathbf{u}^n_h -\mathbf{u}^n = \mathcal{L}_h(\mathbf{u}^n_h - \mathcal{J}_h \tilde{\mathbf{u}}^n) + (\mathcal{L}_h\circ\mathcal{J}_h - I) \tilde{\mathbf{u}}^n
\]
 in which the operators $\mathcal{L}_h$ and $\mathcal{J}_h$ are given by
 \begin{equation*}
 \mathcal{L}_h\mathbf{u}^n_h =\left[
\begin{array}{c}
\mathcal{L}^V_h u^n_h\\
\mathcal{L}^V_h v^n_h
\end{array}
\right],
\quad
\mathcal{J}_h \tilde{\mathbf{u}}^n =\left[
\begin{array}{c}
\mathcal{L}^{V^\ast}_h \tilde{u}^n \\
\mathcal{I}_h \tilde{v}^n
\end{array}
\right].
\end{equation*}
By using the triangle inequality we get
\begin{equation} \label{main_inequality}
\begin{split}
	\|\mathcal{L}_h\mathbf{u}^n_h -\mathbf{u}^n\|_X &\leqslant \|\mathcal{L}_h(\mathbf{u}^n_h - \mathcal{J}_h \tilde{\mathbf{u}}^n)\|_X + \|(\mathcal{L}_h\circ\mathcal{J}_h - I) \tilde{\mathbf{u}}^n\|_X \\
	&\leqslant \max\{C_{\mathrm{H}}, C_{\mathrm{V}}\}\|\mathbf{u}^n_h - \mathcal{J}_h \tilde{\mathbf{u}}^n\|_{X_h} + \|(\mathcal{L}_h\circ\mathcal{J}_h - I) \tilde{\mathbf{u}}^n\|_X.
\end{split}
\end{equation}

\textbf{Step 2}. Due to the following inequality
\begin{equation*}
\small
\begin{split}
	&\|(\mathcal{L}^{V^\ast}_h -\mathcal{I}_h) \tilde{u}^n\|_{V_h} \\
	&= \sup_{\|\omega_h\|_{V_h} =1}\mathrm{p}_h((\mathcal{L}^{V^\ast}_h -\mathcal{I}_h) \tilde{u}^n, \omega_h) \\ 
	&= \sup_{\|\omega_h\|_{V_h} =1}\left\{\mathrm{p}_h(\mathcal{L}^{V^\ast}_h  \tilde{u}^n, \omega_h) -\mathrm{p}_h(\mathcal{I}_h \tilde{u}^n, \omega_h) \right\} \\ 
	&= \sup_{\|\omega_h\|_{V_h} =1}\left\{\mathrm{p}(\tilde{u}^n, \mathcal{L}^{V}_h\omega_h)- \mathrm{p}(\mathcal{L}^{V}_h\circ\mathcal{I}_{h}\tilde{u}^n, \mathcal{L}^{V}_h\omega_h) + \mathrm{p}(\mathcal{L}^{V}_h\circ\mathcal{I}_{h}\tilde{u}^n, \mathcal{L}^{V}_h\omega_h) -\mathrm{p}_h(\mathcal{I}_h \tilde{u}^n, \omega_h) \right\} \\
	&\leqslant \sup_{\|\omega_h\|_{V_h} =1}\left\{\mathrm{p}((I-\mathcal{L}^{V}_h\circ\mathcal{I}_{h})\tilde{u}^n, \mathcal{L}^{V}_h\omega_h) \right\} + \sup_{\|\omega_h\|_{V_h} =1}\left\{\mathrm{p}(\mathcal{L}^{V}_h\circ\mathcal{I}_{h}\tilde{u}^n, \mathcal{L}^{V}_h\omega_h) -\mathrm{p}_h(\mathcal{I}_h \tilde{u}^n, \omega_h) \right\} \\
	&\leqslant C_{\mathrm{V}}\|(I-\mathcal{L}^{V}_h\circ\mathcal{I}_{h})\tilde{u}^n\|_V +  \sup_{\|\omega_h\|_{V_h} =1}\left\{\mathrm{p}(\mathcal{L}^{V}_h\circ\mathcal{I}_{h}\tilde{u}^n, \mathcal{L}^{V}_h\omega_h) -\mathrm{p}_h(\mathcal{I}_h \tilde{u}^n, \omega_h) \right\}\\
	&\leqslant C_{\mathrm{V}}\epsilon^{\mathrm{ip}}_h +\epsilon^{\mathrm{p}}_h,
\end{split}
\end{equation*}
we obtain
\begin{equation}
\label{part_1}
\begin{split}
	 \|(\mathcal{L}_h\circ\mathcal{J}_h - I) \tilde{\mathbf{u}}^n\|_X & = \|(\mathcal{L}^V_h\circ\mathcal{L}^{V^\ast}_h - I) \tilde{u}^n\|_V + \|(\mathcal{L}^V_h\circ\mathcal{I}_h - I) \tilde{v}^n\|_H \\
	 &=  \|(\mathcal{L}^V_h\circ\mathcal{L}^{V^\ast}_h -\mathcal{L}^V_h\circ\mathcal{I}_h + \mathcal{L}^V_h\circ\mathcal{I}_h- I) \tilde{u}^n\|_V + \|(\mathcal{L}^V_h\circ\mathcal{I}_h - I) \tilde{v}^n\|_H\\
	 &\leqslant \|\mathcal{L}^V_h(\mathcal{L}^{V^\ast}_h -\mathcal{I}_h) \tilde{u}^n\|_V + \|(\mathcal{L}^V_h\circ\mathcal{I}_h- I) \tilde{u}^n\|_V + \|(\mathcal{L}^V_h\circ\mathcal{I}_h - I) \tilde{v}^n\|_H \\
	 %&\leqslant C_{\mathrm{V}} {\color{red}\|(\mathcal{L}^{V^\ast}_h -\mathcal{I}_h) \tilde{u}^n\|_{V_h} }+ \|(\mathcal{L}^V_h\circ\mathcal{I}_h- I) \tilde{u}^n\|_V + \|(\mathcal{L}^V_h\circ\mathcal{I}_h - I) \tilde{v}^n\|_H\\
	 &\leqslant C_{\mathrm{V}} \|(\mathcal{L}^{V^\ast}_h -\mathcal{I}_h) \tilde{u}^n\|_{V_h} + \epsilon^{\mathrm{ip}}_h\\
	 %&\leqslant C^2_{\mathrm{V}}\|(I-\mathcal{L}^{V}_h\circ\mathcal{I}_{h})\tilde{u}^n\|_V + C_{\mathrm{V}} \sup_{\|\omega_h\|_{V_h} =1}\left\{\hat{a}(\mathcal{L}^{V}_h\circ\mathcal{I}_{h}\tilde{u}^n, \mathcal{L}^{V}_h\omega_h) -\hat{a}_h(\mathcal{I}_h \tilde{u}^n, \omega_h) \right\}\\
	 %&\qquad + \|(\mathcal{L}^V_h\circ\mathcal{I}_h- I) \tilde{u}^n\|_V + \|(\mathcal{L}^V_h\circ\mathcal{I}_h - I) \tilde{v}^n\|_H\\
	 &\leqslant C^2_{\mathrm{V}}\epsilon^{\mathrm{ip}}_h + C_{\mathrm{V}} \epsilon^{\mathrm{p}}_h + \epsilon^{\mathrm{ip}}_h\\
	&\leqslant (1+C^2_{\mathrm{V}})\left\{ \epsilon^{\mathrm{ip}}_h +\epsilon^{\mathrm{p}}_h \right\}.
\end{split}
\end{equation}

\textbf{Step 3}. Let $\mathbf{e}_h^n := \mathbf{u}^n_h - \mathcal{J}_h \tilde{\mathbf{u}}^n$ denote the discrete error. We will derive the representation of $\mathbf{e}_h^n$.

The fully discrete scheme \eqref{full_discretization} can be rewritten as
\begin{equation} \label{full_discretization_2}
\hat{P}_{+}\mathbf{u}_h^{n+1} = \hat{P}_{-}\mathbf{u}_h^{n} + \frac{\tau}{2}\left(g_h^{n} + g_h^{n+1}\right)+\frac{\tau^2}{4}\left[
\begin{array}{c}
f_h^{n} - f_h^{n+1}\\
-(B_h+\gamma^{n+1})\left(f_h^n - f_h^{n+1}\right)
\end{array}
\right]
\end{equation}
with the discrete operators $\hat{P}_{\pm} = I \pm \frac{\tau}{2}S_h$ and $g_h^{n} = g_h(t_n, \mathbf{u}_h^n)$.
Thanks to $\frac{\tau}{2}\left( \hat{\alpha}\frac{\hat{C}_{\mathrm{em}}}{2} + \hat{\beta} \right) < 1$, the assertions that $\hat{P}_{+}$ is invertible holds true by \cite[Lemma 2.14]{Hipp2017}, \eqref{full_discretization_2} can be reformulated into
\begin{equation}\label{full_discretization_3}
\mathbf{u}_h^{n+1} = \hat{P}\mathbf{u}_h^{n} + \frac{\tau}{2}\hat{P}^{-1}_{+}\left(g_h^{n} + g_h^{n+1}\right)+\frac{\tau^2}{4}\hat{P}^{-1}_{+}\left[
\begin{array}{c}
f_h^{n} - f_h^{n+1}\\
-(B_h+\gamma^{n+1})\left(f_h^n - f_h^{n+1}\right)
\end{array}
\right]
\end{equation}
with $\hat{P}=\hat{P}^{-1}_{+}\hat{P}_{-}$.

Since the exact solution $\tilde{\mathbf{u}}(t)$ satisfies \eqref{eq:evolution}, we have
\begin{equation*}
\begin{aligned}
\frac{\tau}{2}(\dot{\tilde{\mathbf{u}}}^{n+1} + \dot{\tilde{\mathbf{u}}}^{n}) = \frac{\tau}{2}\left( -S(\tilde{\mathbf{u}}^{n+1} + \tilde{\mathbf{u}}^{n}) + \tilde{g}^{n+1}+\tilde{g}^{n}\right)
\end{aligned}
\end{equation*}
with $\tilde{g}^{n} = g( t_n, \tilde{\mathbf{u}}^n)$.
By inserting $\tilde{\mathbf{u}}(t)$ into \eqref{CN} we get
\begin{equation} \label{CN_defect}
	\tilde{\mathbf{u}}^{n+1} = \tilde{\mathbf{u}}^n + \frac{\tau}{2}\left( -S(\tilde{\mathbf{u}}^{n+1} + \tilde{\mathbf{u}}^{n}) + \tilde{g}^{n} + \tilde{g}^{n+1} \right) - \Lambda^{n+1}
\end{equation}
in which $\Lambda^{n+1}$ is a defect.
Then this defect can be given by
\begin{equation*}
	 \Lambda^{n+1} = \frac{\tau}{2}(\dot{\tilde{\mathbf{u}}}^{n+1} + \dot{\tilde{\mathbf{u}}}^{n}) - (\tilde{\mathbf{u}}^{n+1} - \tilde{\mathbf{u}}^n)=\frac{\tau}{2}(\dot{\tilde{\mathbf{u}}}^{n+1} + \dot{\tilde{\mathbf{u}}}^{n})+ \int_{t_n}^{t_{n+1}}\dot{\tilde{\mathbf{u}}}(s)\mathrm{d}s.
\end{equation*}

Taking the operator $\mathcal{L}^{\ast}_h \left[\upsilon,\omega \right]^\top = \left[ \mathcal{L}^{V^\ast}_h \upsilon, \mathcal{L}^{H^\ast}_h\omega \right]^\top$ on both sides of \eqref{CN_defect} yields
\begin{equation*}
	\mathcal{L}^{\ast}_h(\tilde{\mathbf{u}}^{n+1} - \tilde{\mathbf{u}}^n) = \frac{\tau}{2}\left( -\mathcal{L}^{\ast}_h\circ S(\tilde{\mathbf{u}}^{n+1} + \tilde{\mathbf{u}}^{n}) + \mathcal{L}^{\ast}_h\tilde{g}^{n} + \mathcal{L}^{\ast}_h\tilde{g}^{n+1} \right) - \mathcal{L}^{\ast}_h\Lambda^{n+1},
\end{equation*}
from which it deduces
\begin{equation}  \label{CN_defect_2}
\begin{split}
	&\mathcal{J}_h (\tilde{\mathbf{u}}^{n+1} - \tilde{\mathbf{u}}^n) - (\mathcal{J}_h -\mathcal{L}^{\ast}_h)(\tilde{\mathbf{u}}^{n+1} - \tilde{\mathbf{u}}^n) \\
	&=\frac{\tau}{2}\left( -S_h\circ \mathcal{J}_h(\tilde{\mathbf{u}}^{n+1} + \tilde{\mathbf{u}}^{n}) + \tilde{g}_h^{n} + \tilde{g}_h^{n+1} \right) -\frac{\tau}{2}\left( -S_h\circ \mathcal{J}_h(\tilde{\mathbf{u}}^{n+1} + \tilde{\mathbf{u}}^{n}) + \tilde{g}_h^{n} + \tilde{g}_h^{n+1} \right)\\
	&\quad + \frac{\tau}{2}\left( -\mathcal{L}^{\ast}_h\circ S(\tilde{\mathbf{u}}^{n+1} + \tilde{\mathbf{u}}^{n}) + \mathcal{L}^{\ast}_h\tilde{g}^{n} + \mathcal{L}^{\ast}_h\tilde{g}^{n+1} \right) - \mathcal{L}^{\ast}_h\Lambda^{n+1}
\end{split}
\end{equation}
with
\begin{equation*}
\widetilde{g}^n_h = g^n_h(t_n, \mathcal{J}_h\tilde{\mathbf{u}}^n)
=
\left[
\begin{array}{c}
0\\
f_h(\mathcal{L}^{V^\ast}_h \tilde{u}^n) - \gamma(t_n)\mathcal{I}_h \tilde{v}^n
\end{array}
\right] 
=
\left[
\begin{array}{c}
0\\
\tilde{f}^n_h - \gamma^n\mathcal{I}_h \tilde{v}^n
\end{array}
\right].
\end{equation*}
By setting
\begin{equation*}
\small
\begin{split}
	\Pi_h^{n+1} = ( \mathcal{L}^{\ast}_h - \mathcal{J}_h)(\tilde{\mathbf{u}}^{n+1} - \tilde{\mathbf{u}}^n) +\frac{\tau}{2}(\mathcal{L}^{\ast}_h\circ S-S_h\circ \mathcal{J}_h)(\tilde{\mathbf{u}}^{n+1} + \tilde{\mathbf{u}}^{n}) -\frac{\tau}{2}(\mathcal{L}^{\ast}_h\tilde{g}^{n} - \tilde{g}_h^{n} + \mathcal{L}^{\ast}_h\tilde{g}^{n+1} - \tilde{g}_h^{n+1} ),
\end{split}
\end{equation*}
\eqref{CN_defect_2} can be expressed as
\begin{equation*}
	\mathcal{J}_h \tilde{\mathbf{u}}^{n+1} 
	= \mathcal{J}_h \tilde{\mathbf{u}}^{n} + \frac{\tau}{2}\left( -S_h\circ \mathcal{J}_h(\tilde{\mathbf{u}}^{n+1} + \tilde{\mathbf{u}}^{n}) + \tilde{g}_h^{n} + \tilde{g}_h^{n+1} \right) - \Pi_h^{n+1} - \mathcal{L}^{\ast}_h\Lambda^{n+1},
\end{equation*}
which is equivalent to the following representation
\begin{equation} \label{CN_defect_3}
	\mathcal{J}_h \tilde{\mathbf{u}}^{n+1} 
	=\hat{P}\circ \mathcal{J}_h \tilde{\mathbf{u}}^{n} + \frac{\tau}{2}\hat{P}^{-1}_{+}\left(\tilde{g}_h^{n} + \tilde{g}_h^{n+1}\right) - \hat{P}^{-1}_{+}\Pi_h^{n+1} - \hat{P}^{-1}_{+}\mathcal{L}^{\ast}_h\Lambda^{n+1}.
\end{equation}

Plugging $\mathcal{J}_h \tilde{\mathbf{u}}$ into \eqref{full_discretization_3} yields
\begin{equation}\label{IMEX_error}
\mathcal{J}_h \tilde{\mathbf{u}}^{n+1} = \hat{P}\circ\mathcal{J}_h\tilde{\mathbf{u}}^{n} + \frac{\tau}{2}\hat{P}^{-1}_{+}\left(\tilde{g}_h^{n} + \tilde{g}_h^{n+1}\right) +
\frac{\tau^2}{4}\hat{P}^{-1}_{+}\left[
\begin{array}{c}
\tilde{f}_h^{n} - \tilde{f}_h^{n+1}\\
-(B_h + \gamma^{n+1}) \left(\tilde{f}_h^n - \tilde{f}_h^{n+1}\right)
\end{array}
\right]
-\Delta^{n+1}_h
\end{equation}
with another defect $\Delta^{n+1}_h$. By comparing to \eqref{CN_defect_3}, we get the 
\begin{equation*}
\Delta^{n+1}_h = \hat{P}^{-1}_{+}\Pi_h^{n+1} + \hat{P}^{-1}_{+}\mathcal{L}^{\ast}_h\Lambda^{n+1} + \Theta_h^{n+1},
\quad
\Theta_h^{n+1} = 
\frac{\tau^2}{4}\hat{P}^{-1}_{+}\left[
\begin{array}{c}
\tilde{f}_h^{n} - \tilde{f}_h^{n+1}\\
-(B_h + \gamma^{n+1}) \left(\tilde{f}_h^n - \tilde{f}_h^{n+1}\right)
\end{array}
\right].
\end{equation*}

By subtracting \eqref{IMEX_error} from \eqref{full_discretization_3} we obtain the following recursion
\begin{equation*}
\small
\mathbf{e}^{n+1}_h = \hat{P}\mathbf{e}^{n}_h + \frac{\tau}{2}\hat{P}^{-1}_{+}\left(g_h^{n} - \tilde{g}_h^{n} + g_h^{n+1} - \tilde{g}_h^{n+1}\right)+
\frac{\tau^2}{4}\hat{P}^{-1}_{+}\left[
\begin{array}{c}
f_h^{n} - \tilde{f}_h^{n} - f_h^{n+1} + \tilde{f}_h^{n+1}\\
-(B_h + \gamma^{n+1}) \left(f_h^{n} - \tilde{f}_h^{n} - f_h^{n+1} + \tilde{f}_h^{n+1}\right)
\end{array}
\right]
-\Delta^{n+1}_h,
\end{equation*}
which implies
\begin{equation} \label{error_recursion}
\begin{split}
\mathbf{e}^{n}_h &= \hat{P}^n\mathbf{e}^{0}_h + \sum_{m=1}^{n}\hat{P}^{n-m}\bigg\{\frac{\tau}{2}\hat{P}^{-1}_{+}\left(g_h^{m} - \tilde{g}_h^{m} + g_h^{m +1} - \tilde{g}_h^{m +1}\right)\\
& \quad+\frac{\tau^2}{4}\hat{P}^{-1}_{+}\left[
\begin{array}{c}
f_h^{m-1} - \tilde{f}_h^{m-1} - f_h^{m} + \tilde{f}_h^{m}\\
-(B_h + \gamma^{m+1}) \left(f_h^{m-1} - \tilde{f}_h^{m-1} - f_h^{m} + \tilde{f}_h^{m}\right)
\end{array}
\right]
-\Delta^{m}_h \bigg\}.
\end{split}
\end{equation}

\textbf{Step 4}. To bound $\|\mathbf{e}_h^n\|_{X_h}$, we need some preparations. By \cite[Lemma 2.14]{Hipp2017} we know
\[
	\|\hat{P}\upsilon\|_{X_{h}}\leqslant e^{\tau \hat{c}_{0}}\|\upsilon\|_{X_{h}}, \quad \|\hat{P}^{-1}_{+}\upsilon\|_{X_{h}}\leqslant \|\upsilon\|_{X_{h}},\quad \forall \upsilon\in X_h
\]
with $\hat{c}_{0} = \hat{ \alpha}\frac{\hat{C}_{\mathrm{em}}}{2} + \hat{\beta}$. Then we have
\begin{equation} \label{step4_1}
\begin{split}
&\left\| \hat{P}^{n-m}\left\{\frac{\tau}{2}\hat{P}^{-1}_{+}\left(g_h^{m} - \tilde{g}_h^{m} + g_h^{m -1} - \tilde{g}_h^{m -1}\right)\right\} \right\|_{X_h} \\
&\leqslant \frac{\tau}{2}e^{(n-m)\tau \hat{c}_{0} }\left( \|\hat{P}^{-1}_{+}\left(g_h^{m} - \tilde{g}_h^{m}\right) \|_{X_h} + \|\hat{P}^{-1}_{+}\left(g_h^{m-1} - \tilde{g}_h^{m-1}\right) \|_{X_h}\right)\\
&\leqslant \frac{\tau}{2}e^{(n-m)\tau \hat{c}_{0} }\left( \|g_h^{m} - \tilde{g}_h^{m} \|_{X_h} + \| g_h^{m-1} - \tilde{g}_h^{m-1} \|_{X_h}\right).
\end{split}
\end{equation}
We note that $\|u^m_h\|_{V_h}  \leqslant\delta$ and 
\begin{equation*}
\begin{split}
	\|\mathcal{L}^{V^\ast}_h \tilde{u}^m\|_{V_h} &= \sup_{\|\omega_h\|_{V_h} =1}\mathrm{p}_h(\mathcal{L}^{V^\ast}_h  \tilde{u}^m, \omega_h) = \sup_{\|\omega_h\|_{V_h} =1}\mathrm{p}(\tilde{u}^m, \mathcal{L}^{V}_h  \omega_h)\leqslant C_{\mathrm{V}}\|\tilde{u}^m\|_V \leqslant\delta.
\end{split}
\end{equation*}
Then by using the Lipschitz-continuity of the discrete nonlinearity $f_h$ we have 
\begin{equation} \label{step4_2}
\begin{split}
\|g_h^{m} - \tilde{g}_h^{m} \|_{X_h}
&=
\left\|\left[
\begin{array}{c}
0\\
f_h(u^m_h) - \gamma(t_m)v^m_h -f_h(\mathcal{L}^{V^\ast}_h \tilde{u}^m) + \gamma(t_m)\mathcal{I}_h \tilde{v}^m
\end{array}
\right]\right\|_{X_h}\\
&\leqslant \left\| f_h(u^m_h) - f_h(\mathcal{L}^{V^\ast}_h \tilde{u}^m) \right\|_{H_h} + \gamma(t_m)\left\| v^m_h - \mathcal{I}_h \tilde{v}^m \right\|_{H_h}\\
&\leqslant \hat{C}_{\delta}\| u^m_h - \mathcal{L}^{V^\ast}_h \tilde{u}^m \|_{V_h} + (\gamma_0+ T)\| v^m_h - \mathcal{I}_h \tilde{v}^m \|_{H_h}\\
&\leqslant (\hat{C}_{\delta} + C_{\gamma})\| \mathbf{e}^m_h \|_{X_h}.
\end{split}
\end{equation}
Combing \eqref{step4_1} and \eqref{step4_2} gives
\begin{equation} \label{step4_3}
\begin{split}
\left\| \hat{P}^{n-m}\left\{\frac{\tau}{2}\hat{P}^{-1}_{+}\left(g_h^{m} - \tilde{g}_h^{m} + g_h^{m -1} - \tilde{g}_h^{m -1}\right)\right\} \right\|_{X_h} \leqslant \frac{\tau}{2}(\hat{C}_{\delta} + C_{\gamma})e^{(n-m)\tau \hat{c}_{0} }\left(\| \mathbf{e}^m_h \|_{X_h} + \| \mathbf{e}^{m-1}_h \|_{X_h} \right).
\end{split}
\end{equation}

Let $\hat{Q}_{\pm} = I \pm \frac{\tau^2}{4}A_h \pm \frac{\tau}{2} B_h$.
For $\frac{\tau^2}{2}\hat{\alpha} + \tau\hat{\beta} < 1$ by \cite[Lemma 4.6]{Leibold2021} we obtain $\hat{Q}_{+}$ is invertible and its inverse $\hat{Q}_{+}^{-1}$ satisfies
\[
	\|\hat{Q}_{+}^{-1}\upsilon\|_{V_h} \leqslant \frac{\sqrt{2}}{\tau}\| \upsilon \|_{H_h}, \quad \|(-I + \hat{Q}_{+}^{-1})\upsilon\|_{H_h} \leqslant \| \upsilon \|_{H_h}, \quad \|(-I+2\hat{Q}_{+}^{-1})\upsilon\|_{V_h} \leqslant e^{\frac{\tau^2}{2}\hat{\alpha} + \tau\hat{\beta} }\| \upsilon \|_{H_h}
\]
for any $\upsilon\in H_h$. Furthermore, we can derive the representations of $\hat{P}^{-1}_{+}$ as 
\[
\hat{P}^{-1}_{+} = 
\left[
\begin{array}{cc}
\hat{Q}_{+}^{-1}(I+\frac{\tau}{2}B_h) & \frac{\tau}{2}\hat{Q}_{+}^{-1}\\
-\frac{2}{\tau}+ \frac{2}{\tau}\hat{Q}_{+}^{-1}(I+\frac{\tau}{2}B_h)  & \hat{Q}_{+}^{-1}
\end{array}
\right] 
\]
and
\begin{equation*}
\begin{split}
\hat{P}&= \hat{P}^{-1}_{+}\hat{P}_{-} =\hat{P}^{-1}_{+}(I-\frac{\tau}{2}S_h) =
\left[
\begin{array}{cc}
-I + \hat{Q}_{+}^{-1}(2I+B_h) & \tau\hat{Q}_{+}^{-1}\\
-\frac{4}{\tau}+ \frac{2}{\tau}\hat{Q}_{+}^{-1}(2I+\tau B_h)  & -I +2\hat{Q}_{+}^{-1}
\end{array}
\right]. 
\end{split}
\end{equation*}
Thus, we have the estimates
\begin{equation}\label{step4_4}
\begin{split}
\left\| \hat{P}^{-1}_{+}\left[
\begin{array}{c}
f_h^{m} - \tilde{f}_h^{m}\\
-B_h\left(f_h^{m} - \tilde{f}_h^{m} \right)
\end{array}
\right] \right\|_{X_h} 
&= 
\left\| 
\left[
\begin{array}{c}
\hat{Q}_{+}^{-1}(f_h^{m} - \tilde{f}_h^{m})\\
(-\frac{2}{\tau} + \frac{2}{\tau}\hat{Q}_{+}^{-1})\left(f_h^{m} - \tilde{f}_h^{m} \right)
\end{array}
\right] 
\right\|_{X_h} 
\\
&\leqslant  \frac{\sqrt{2}}{\tau}\|f_h^{m} - \tilde{f}_h^{m}\|_{H_h} + \frac{2}{\tau}\|f_h^{m} - \tilde{f}_h^{m}\|_{H_h}\\
&\leqslant  \frac{2\sqrt{3}}{\tau}\hat{C}_{\delta}\|\mathbf{e}^{m}_h\|_{X_h}
\end{split}
\end{equation}
and 
\begin{equation} \label{step4_5}
\begin{split}
\left\| \hat{P}^{-1}_{+}\left[
\begin{array}{c}
0\\
-\gamma^{m+1} \left(f_h^{m} - \tilde{f}_h^{m} \right)
\end{array}
\right] \right\|_{X_h}
&= 
\left\| 
\gamma^{m+1}\left[
\begin{array}{c}
\frac{\tau}{2}\hat{Q}_{+}^{-1}(f_h^{m} - \tilde{f}_h^{m})\\
\hat{Q}_{+}^{-1}\left(f_h^{m} - \tilde{f}_h^{m} \right)
\end{array}
\right] 
\right\|_{X_h} 
\\
&\leqslant C_{\gamma}(\frac{\sqrt{2}}{2} + \frac{\sqrt{2}}{\tau})\hat{C}_{\delta}\|\mathbf{e}^{m}_h\|_{X_h}.
\end{split}
\end{equation}

\textbf{Step 5}. Taking the norm on \eqref{error_recursion}, using the triangle inequality, \eqref{step4_3}, \eqref{step4_4} and \eqref{step4_5} gives
\begin{equation*}
\begin{split}
\left\| \mathbf{e}^{n}_h \right\|_{X_h}  &\leqslant  \left\| \hat{P}^n\mathbf{e}^{0}_h +  \sum_{m=1}^{n}\hat{P}^{n-m}\Delta^{m}_h\right\|_{X_h} + \sum_{m=1}^{n} \left\|\hat{P}^{n-m}\bigg\{\frac{\tau}{2}\hat{P}^{-1}_{+}\left(g_h^{m} - \tilde{g}_h^{m} + g_h^{m +1} - \tilde{g}_h^{m +1}\right) \bigg\} \right\|_{X_h}\\
& \quad+\frac{\tau^2}{4}\sum_{m=1}^{n} \left\|\hat{P}^{n-m}\circ\hat{P}^{-1}_{+}\left[
\begin{array}{c}
f_h^{m-1} - \tilde{f}_h^{m-1} - f_h^{m} + \tilde{f}_h^{m}\\
-B_h \left(f_h^{m-1} - \tilde{f}_h^{m-1} - f_h^{m} + \tilde{f}_h^{m}\right)
\end{array}
\right] \right\|_{X_h}\\
& \quad+\frac{\tau^2}{4}\sum_{m=1}^{n} \left\|\hat{P}^{n-m}\circ\hat{P}^{-1}_{+}\left[
\begin{array}{c}
0\\
- \gamma^m \left(f_h^{m-1} - \tilde{f}_h^{m-1} - f_h^{m} + \tilde{f}_h^{m}\right)
\end{array}
\right] \right\|_{X_h}\\
&\leqslant \tau \hat{C}_{re}\sum_{m=1}^{n}e^{(n-m)\tau \hat{c}_{0} }\|\mathbf{e}^{m}_h\|_{X_h} + \left\| \hat{P}^n\mathbf{e}^{0}_h +  \sum_{m=1}^{n}\hat{P}^{n-m}\Delta^{m}_h\right\|_{X_h}
\end{split}
\end{equation*}
with $\hat{C}_{re} = (1+\sqrt{3})\hat{C}_{\delta} +C_{\gamma} + \frac{3\sqrt{2}}{4}\hat{C}_{\delta}C_{\gamma}$.
Thanks to $\tau \hat{C}_{\mathrm{re}} < 1$, we obtain
\begin{equation} \label{recursion_norm}
\begin{split}
\left\| \mathbf{e}^{n}_h \right\|_{X_h}  &\leqslant  e^{\frac{n\tau \hat{C}_{re}}{1- \tau \hat{C}_{\mathrm{re}}}}\left\| \hat{P}^n\mathbf{e}^{0}_h +  \sum_{m=1}^{n}\hat{P}^{n-m}\Delta^{m}_h\right\|_{X_h} \leqslant e^{\frac{n\tau \hat{C}_{re}}{1- \tau \hat{C}_{\mathrm{re}}}}\left(\left\| \hat{P}^n\mathbf{e}^{0}_h \right\|_{X_h} + \left\|  \sum_{m=1}^{n}\hat{P}^{n-m}\Delta^{m}_h\right\|_{X_h} \right).
\end{split}
\end{equation}
by utilizing a discrete Gronwall inequality and the triangle inequality. 

Due to
\begin{equation*}
\small
\begin{split}
\left\|\Lambda^{m}\right\|_{X_h} = \left\|\tau^3\int_{0}^{1}\frac{s(1-s)}{2} \frac{\mathrm{d}^3}{\mathrm{d}t^3} \mathbf{u}(t_m + \tau s)\mathrm{d}s\right\|_{X_h} \leqslant \tau^3\left(  \left\|\frac{\mathrm{d}^4 u}{\mathrm{d}t^4}\right\|_{L^{\infty}([t_{m-1}, t_m]; H)} + \left\|\frac{\mathrm{d}^3 u}{\mathrm{d}t^3}\right\|_{L^{\infty}([t_{m-1}, t_m]; V)} \right)
\end{split}
\end{equation*}
and
\begin{equation*} 
\small
\begin{split}
\left\|\Pi_h^{m}\right\|_{X_h} &= \|(\mathcal{L}^{V^\ast}_h -\mathcal{I}_h) (\tilde{v}^{m+1}- \tilde{v}^{m})\|_{H_h} + \frac{\tau}{2}\|(\mathcal{L}^{V^\ast}_h -\mathcal{I}_h) (\tilde{v}^{m+1} + \tilde{v}^{m})\|_{V_h} \\
&\quad + \frac{\tau}{2}\left\| \mathcal{L}^{H^\ast}_h\left(A (\tilde{u}^{m+1} + \tilde{u}^{m}) + B(\tilde{v}^{m+1} + \tilde{v}^{m}) \right) - \left(A_h \mathcal{L}^{H^\ast}_h(\tilde{u}^{m+1} + \tilde{u}^{m}) + B_h(\tilde{v}^{m+1} + \tilde{v}^{m})\mathcal{I}_h\right) \right\|_{H_h} \\
&\quad + \frac{\tau}{2}\left\|\mathcal{L}^{H^\ast}_h f(\tilde{u}^{m+1}) - f_h(\mathcal{L}^{V^\ast}_h \tilde{u}^{m+1}) \right\|_{H_h} + \frac{\tau}{2}\left\|\mathcal{L}^{H^\ast}_h f(\tilde{u}^{m}) - f_h(\mathcal{L}^{V^\ast}_h \tilde{u}^{m}) \right\|_{H_h}\\
&\quad+ \gamma^{m}\frac{\tau}{2}\|(\mathcal{L}^{H^\ast}_h -\mathcal{I}_h)(\tilde{v}^{m+1} + \tilde{v}^{m})  \|_{H_h}\\
&\leqslant \tau(C_{\mathrm{V}}\epsilon^{\mathrm{ip}}_h + \epsilon^{\mathrm{p}}_h) +\tau (C_{\mathrm{V}}\epsilon^{\mathrm{ip}}_h + \epsilon^{\mathrm{p}}_h) + C\tau(\epsilon^{\mathrm{ip}}_h+\epsilon^{\mathrm{p}}_h + \epsilon^{\mathrm{b}}_h) + \tau \epsilon^{\mathrm{data}}_h ++\tau C_\gamma(C_{\mathrm{H}}\epsilon^{\mathrm{ip}}_h + \epsilon^{\mathrm{p}}_h)\\
&\leqslant C_1 \tau\epsilon^{\mathrm{total}}_h,
\end{split}
\end{equation*}
we have
\begin{equation} \label{defect_estimate_1}
\begin{split}
\left\|\sum_{m=1}^{n}\hat{P}^{n-m}\Delta^{m}_h \right\|_{X_h} &= \left\|\sum_{m=1}^{n}\hat{P}^{n-m}\left( \hat{P}^{-1}_{+}\Pi_h^{m} + \hat{P}^{-1}_{+}\mathcal{L}^{\ast}_h\Lambda^{m} + \Theta_h^{m} \right)\right\|_{X_h} \\
&\leqslant \sum_{m=1}^{n}e^{(n-m)\tau \hat{c}_{0} }\left(\left\|\Pi_h^{m}\right\|_{X_h} + \max\{C_{\mathrm{H}}, C_{\mathrm{V}}\}\left\|\Lambda^{m}\right\|_{X_h}\right) + \left\|\sum_{m=1}^{n}\hat{P}^{n-m}\Theta_h^{m} \right\|_{X_h}\\
&\leqslant C_2 e^{n\tau \hat{c}_{0} }\left(\tau^2 + \epsilon^{\mathrm{total}}_h \right) + \left\|\sum_{m=1}^{n}\hat{P}^{n-m}\Theta_h^{m} \right\|_{X_h}.
\end{split}
\end{equation}

It remains to bound the second term on the left hand side of \eqref{defect_estimate_1}. We set $\Theta_h^{m} = \Phi_h^{m} + \Psi_h^{m}$ with
\begin{equation*}
\begin{split}
\Phi_h^{m} 
= 
\frac{\tau^2}{4}\hat{P}^{-1}_{+}\left[
\begin{array}{c}
\tilde{f}_h^{m-1} - \tilde{f}_h^{m}\\
-B_h \left(\tilde{f}_h^{m-1} - \tilde{f}_h^{m}\right)
\end{array}
\right], \quad
\Psi_h^{m} 
= 
\frac{\tau^2}{4}\hat{P}^{-1}_{+}\left[
\begin{array}{c}
0\\
-\gamma^{m} \left(\tilde{f}_h^{m-1} - \tilde{f}_h^{m}\right)
\end{array}
\right].
\end{split}
\end{equation*}
We note that
\begin{equation*}
\small
\begin{split}
\hat{Q}_{+}^{-1}\left(\tilde{f}_h^{m-1} - \tilde{f}_h^{m} \right) = \hat{Q}_{+}^{-1}\left(\tilde{f}_h^{m-1} - \mathcal{L}^{H^\ast}_h\tilde{f}^{m-1} - \tilde{f}_h^{m} +  \mathcal{L}^{H^\ast}_h\tilde{f}^{m}\right) +\hat{Q}_{+}^{-1}\circ \mathcal{L}^{H^\ast}_h\left(\tilde{f}^{m-1} - \tilde{f}^{m}\right)
\end{split}
\end{equation*}
and
\begin{equation*}
\small
\begin{split}
&\left(-2 + 2\hat{Q}_{+}^{-1}\right)\left(- \mathcal{L}^{H^\ast}_h\tilde{f}^{m-1} +  \mathcal{L}^{H^\ast}_h\tilde{f}^{m}\right) + \left(-I+2\hat{Q}_{+}^{-1}\right)\circ \mathcal{L}^{H^\ast}_h\left(\tilde{f}^{m-1} - \tilde{f}^{m}\right) - \mathcal{L}^{H^\ast}_h\left(\tilde{f}^{m-1} - \tilde{f}^{m}\right)\\
&=-2\hat{Q}_{+}^{-1} \circ\mathcal{L}^{H^\ast}_h\left(\tilde{f}^{m-1} -\tilde{f}^{m}\right)+\left(-I+2\hat{Q}_{+}^{-1}\right)\circ \mathcal{L}^{H^\ast}_h\left(\tilde{f}^{m-1} - \tilde{f}^{m}\right)+\mathcal{L}^{H^\ast}_h\left(\tilde{f}^{m-1} - \tilde{f}^{m}\right)\\
&=\left(-2\hat{Q}_{+}^{-1} + (-I+2\hat{Q}_{+}^{-1})\right)\mathcal{L}^{H^\ast}_h\left(\tilde{f}^{m-1} - \tilde{f}^{m}\right) +\mathcal{L}^{H^\ast}_h\left(\tilde{f}^{m-1} - \tilde{f}^{m}\right)\\
&=0.
\end{split}
\end{equation*}
Then we have the following decomposition
\begin{equation} \label{decomposition_1}
\sum_{m=1}^{n}\hat{P}^{n-m}\Phi_h^{m}
=
\frac{\tau^2}{4}
\sum_{m=1}^{n}\hat{P}^{n-m}\left[
\begin{array}{c}
\hat{Q}_{+}^{-1}\left(\tilde{f}_h^{m-1} - \tilde{f}_h^{m}\right)\\
(-\frac{2}{\tau} + \frac{2}{\tau}\hat{Q}_{+}^{-1})\left(\tilde{f}_h^{m-1} - \tilde{f}_h^{m}\right)
\end{array}
\right] = \sum_{m=1}^{n}\hat{P}^{n-m}\left(\Xi^{m}_{1} + \Xi^{m}_2 + \Xi^{m}_3 \right)
\end{equation} 
with
\[
\Xi^m_1 = \frac{\tau^2}{4}
\left[
\begin{array}{c}
\hat{Q}_{+}^{-1}\left(\tilde{f}_h^{m-1} - \mathcal{L}^{H^\ast}_h\tilde{f}^{m-1} - \tilde{f}_h^{m} +  \mathcal{L}^{H^\ast}_h\tilde{f}^{m}\right)\\
(-\frac{2}{\tau} + \frac{2}{\tau}\hat{Q}_{+}^{-1})\left(\tilde{f}_h^{m-1} - \mathcal{L}^{H^\ast}_h \tilde{f}^{m-1} - \tilde{f}_h^{m} +  \mathcal{L}^{H^\ast}_h\tilde{f}^{m}\right)
\end{array}
\right]
\]
and
\[
\Xi^m_2 = \frac{\tau}{4}
\left[
\begin{array}{c}
\tau\hat{Q}_{+}^{-1}\circ \mathcal{L}^{H^\ast}_h\left(\tilde{f}^{m-1} - \tilde{f}^{m}\right)\\
(-I+2\hat{Q}_{+}^{-1})\circ \mathcal{L}^{H^\ast}_h\left(\tilde{f}^{m-1} - \tilde{f}^{m}\right)
\end{array}
\right], 
\quad
\Xi^m_3 =\frac{\tau}{4}
\left[
\begin{array}{c}
0\\
-\mathcal{L}^{H^\ast}_h\left(\tilde{f}^{m-1} - \tilde{f}^{m}\right)
\end{array}
\right]. 
\]
Since we have
\begin{equation*}
\begin{split}
\hat{P}\Xi^{m-1}_3
=
\left[
\begin{array}{c}
-\tau \hat{Q}_{+}^{-1}\mathcal{L}^{H^\ast}_h\left(\tilde{f}^{m-2} - \tilde{f}^{m-1}\right)\\
(I-2\hat{Q}_{+}^{-1})\circ\mathcal{L}^{H^\ast}_h\left(\tilde{f}^{m-2} - \tilde{f}^{m-1}\right)
\end{array}
\right] = -\Xi^{m-1}_2, 
\end{split}
\end{equation*}
the decomposition \eqref{decomposition_1} can be reformulated as
\begin{equation} \label{decomposition_2}
\begin{split}
\sum_{m=1}^{n}\hat{P}^{n-m}\Phi_h^{m}
&= \sum_{m=1}^{n}\hat{P}^{n-m}\Xi^{m}_{1} + \hat{P}^{n-1}\Xi^{1}_2 + \sum_{m=2}^{n}\hat{P}^{n-m}\Xi^{m}_{2}  + \sum_{m=1}^{n-1}\hat{P}^{n-m}\Xi^{m}_3 + \Xi^{n}_3\\
&= \sum_{m=1}^{n}\hat{P}^{n-m}\Xi^{m}_{1} + \hat{P}^{n-1}\Xi^{1}_2 + \sum_{m=2}^{n}\hat{P}^{n-m}\Xi^{m}_{2}  + \sum_{m=2}^{n}\hat{P}^{n-m}(\hat{P}\Xi^{m-1}_3) + \Xi^{n}_3\\
&= \sum_{m=1}^{n}\hat{P}^{n-m}\Xi^{m}_{1} + \hat{P}^{n-1}\Xi^{1}_2 + \sum_{m=2}^{n}\hat{P}^{n-m}(\Xi^{m}_{2} - \Xi^{m-1}_2) + \Xi^{n}_3.
\end{split}
\end{equation} 
By the properties of $\hat{Q}_{+}^{-1}$ and $\mathcal{L}^{H^\ast}_h$, we have
\[
\left\|\Xi^{m}_{1}\right\|_{X_h}
\leqslant C_3\tau \epsilon^{\mathrm{data}}_h,
\quad
\left\|\Xi^{1}_{2}\right\|_{X_h}
\leqslant C_4\tau^2,
\quad
\left\|\Xi^{n}_{3}\right\|_{X_h}
\leqslant C_5\tau^2
\]
and
\begin{equation*}
\small
\begin{split}
\left\|\Xi^{m}_{2} - \Xi^{m-1}_2\right\|_{X_h} &=\frac{\tau}{4}
\left\|\left[
\begin{array}{c}
\tau\hat{Q}_{+}^{-1}\circ \mathcal{L}^{H^\ast}_h\left(-\tilde{f}^{m-2} + 2\tilde{f}^{m-1} - \tilde{f}^{m}\right)\\
(-I+2\hat{Q}_{+}^{-1})\circ \mathcal{L}^{H^\ast}_h\left(-\tilde{f}^{m-2} + 2\tilde{f}^{m-1} - \tilde{f}^{m}\right)
\end{array}
\right] \right\|_{X_h}\\
&\leqslant C_6\tau \left\| -\tilde{f}^{m-2} + 2\tilde{f}^{m-1} - \tilde{f}^{m} \right\|_{H}\\
&\leqslant C_7\tau^3\left\|\frac{\mathrm{d}^2}{\mathrm{d}t^2} f(u(t))\right\|_{L^{\infty}([t_{m-1}, t_m]; H)}\\
&= C_7\tau^3\left\|\frac{\mathrm{d}^4}{\mathrm{d}t^4} u(t) + A\left(\frac{\mathrm{d}^2}{\mathrm{d}t^2} u(t)\right)+ B\left(\frac{\mathrm{d}^3}{\mathrm{d}t^3} u(t)\right) +\frac{\mathrm{d}^2}{\mathrm{d}t^2} \left(\gamma(t)u(t)\right)\right\|_{L^{\infty}([t_{m-1}, t_m]; H)}\\
&\leqslant (C_8+C_{\gamma})\tau^3\left(  \left\|\frac{\mathrm{d}^4 u}{\mathrm{d}t^4}\right\|_{L^{\infty}([t_{m-1}, t_m]; H)} + \left\|\frac{\mathrm{d}^3 u}{\mathrm{d}t^3}\right\|_{L^{\infty}([t_{m-1}, t_m]; H)} + \left\|\ddot{u}\right\|_{L^{\infty}([t_{m-1}, t_m]; H)} + \left\| \ddot{u}\right\|_{L^{\infty}([t_{m-1}, t_m]; H)}\right).
\end{split}
\end{equation*}
Taking the norm and inserting all the corresponding bounds into \eqref{decomposition_2} give
\begin{equation} \label{defect_estimate_2}
\left\|\sum_{m=1}^{n}\hat{P}^{n-m}\Phi_h^{m} \right\|_{X_h} \leqslant C_9e^{n\tau \hat{c}_{0} }\left( \epsilon^{\mathrm{data}}_h + \tau^2\right).
\end{equation}

By the following estimate
\begin{equation*}
\begin{split}
\left\|\Psi_h^{m} \right\|_{X_h}
&=\frac{\tau^2}{4}\gamma^{m}\left\| 
\left[
\begin{array}{c}
\frac{\tau}{2}\hat{Q}_{+}^{-1}\left(\tilde{f}_h^{m-1} - \tilde{f}_h^{m}\right)\\
\hat{Q}_{+}^{-1}\left(\tilde{f}_h^{m-1} - \tilde{f}_h^{m}\right)
\end{array}
\right] 
\right\|_{X_h}\\ 
&\leqslant \frac{\tau^2}{4}(\frac{\sqrt{2}}{2} + \frac{\sqrt{2}}{\tau})C_{\gamma}\left\|\tilde{f}_h^{m-1} - \tilde{f}_h^{m}\right\|_{H_h}\\
&=\frac{\tau^2}{4}(\frac{\sqrt{2}}{2} + \frac{\sqrt{2}}{\tau})C_{\gamma}\left\|f_h(\mathcal{L}^{V^\ast}_h \tilde{u}^{m-1}) -f_h(\mathcal{L}^{V^\ast}_h \tilde{u}^m)\right\|_{H_h}\\
&\leqslant\frac{\tau^2}{4}(\frac{\sqrt{2}}{2} + \frac{\sqrt{2}}{\tau})C_{\gamma}\hat{C}_{\delta}C_{\mathrm{V}}\left\|u(t_{m-1}) -u(t_{m})\right\|_{V}\\
&\leqslant\frac{\tau^2}{4}(\frac{\sqrt{2}}{2}\tau + \sqrt{2})C_{\gamma}\hat{C}_{\delta}C_{\mathrm{V}}\left\|\dot{u}\right\|_{L^{\infty}([t_{m-1}, t_m]; V)} \\
&\leqslant\frac{3\sqrt{2}\tau^2}{8}C_{\gamma}\hat{C}_{\delta}C_{\mathrm{V}}\left\|\dot{u}\right\|_{L^{\infty}([t_{m-1}, t_m]; V)}
\end{split}
\end{equation*}
and setting $C_{10} = \frac{3\sqrt{2}}{8}\hat{C}_{\delta}C_{\mathrm{V}}\left\|\dot{u}\right\|_{L^{\infty}([t_{m-1}, t_m]; V)}$, we get
\begin{equation} \label{defect_estimate_3}
\left\|\sum_{m=1}^{n}\hat{P}^{n-m}\Psi_h^{m} \right\|_{X_h} \leqslant C_{10}C_{\gamma}e^{n\tau \hat{c}_{0} }\tau.
\end{equation}

Inserting \eqref{defect_estimate_1}, \eqref{defect_estimate_2}, \eqref{defect_estimate_3} and 
\begin{equation*} 
\left\| \hat{P}^n\mathbf{e}^{0}_h \right\|_{X_h} \leqslant  e^{n\tau \hat{c}_{0} } \left\|\mathbf{e}^{0}_h \right\|_{X_h}  \leqslant  e^{n\tau \hat{c}_{0} }\epsilon^{\mathrm{data}}_h
\end{equation*}
into \eqref{recursion_norm} yields
\begin{equation} \label{part_2}
\begin{split}
\left\| \mathbf{e}^{n}_h \right\|_{X_h} &\leqslant \exp\left\{\frac{n\tau \hat{C}_{re}}{1- \tau \hat{C}_{\mathrm{re}}}\right\}\left(e^{n\tau \hat{c}_{0} }\epsilon^{\mathrm{data}}_h+ C_2 e^{n\tau \hat{c}_{0} }\left(\tau^2 + \epsilon^{\mathrm{total}}_h \right) + C_9 e^{n\tau \hat{c}_{0} }\left( \epsilon^{\mathrm{data}}_h + \tau^2\right) +C_{10}C_{\gamma}e^{n\tau \hat{c}_{0} }\tau \right)\\
&=C_{11}e^{M n\tau}\left( \epsilon^{\mathrm{total}}_h + \tau^2 + C_{\gamma}\tau \right)
\end{split}
\end{equation}
with $M = \frac{ \hat{C}_{re}}{1- \tau \hat{C}_{\mathrm{re}}} + \hat{c}_{0}$.

\textbf{Step 6}. By \eqref{part_1} and \eqref{part_2}, it deduce from \eqref{main_inequality} that
\[
	\|\mathcal{L}_h\mathbf{u}^n_h -\mathbf{u}^n\|_X \leqslant Ce^{M n\tau}\left( \epsilon^{\mathrm{total}}_h + \tau^2 + C_{\gamma}\tau \right).
\]
Therefore, we complete the proof.
\end{proof}

%%%%%%%%%%%%%%%%%%%%%%%%%%%%%%%%%%%%%%%%%%%%%%%%%%%%%%%%%%%%%%%%%%%%%%%%%%%%%%%%
\section{Conclusions}
\label{sed:conclusion}
In this paper, we have found from our theoretical results that one actually cannot expect $o(\tau^2)$ convergence due to the effect of the nonautonomous damping.
In numerical simulations, the reduced convergence rate can also be observed. Moreover, we have developed a revised IMEX scheme for the semilinear wave-type equations with time-varying dampings. 
Numerical tests are presented to substantiate the validity and efficiency of the revised numerical method. Notably, the numerical results confirms that the second order $o(\tau^2)$ of the error can be reached.

%%%%%%%%%%%%%%%%%%%%%%%%%%%%%%%%%%%%%%%%%%%%%%%%%%%%%%%%%%%%%%%%%%%%%%%%%%%%%%%%

\appendix

\section{Crank-Nicolson scheme for Equation~\eqref{eq:evolution}}
\label{sec:CN}
The Crank-Nicolson scheme for \eqref{eq:evolution} is of the form
\begin{equation}\label{CN}
	\mathbf{u}^{n+1} = \mathbf{u}^n + \frac{\tau}{2}\left( -S(\mathbf{u}^{n} + \mathbf{u}^{n+1}) + g^{n} + g^{n+1}\right),
\end{equation}
which can be written as
\begin{equation}\label{CN_2}
	P_{+}\mathbf{u}^{n+1} = P_{-}\mathbf{u}^{n} + \frac{\tau}{2}\left(g^{n} + g^{n+1}\right)
\end{equation}
with the operator $P_{\pm} = I\pm  \frac{\tau}{2}S$ and $g^{n} = g(t_{n}, \mathbf{u}^n)$.

The two components of \eqref{CN} have the form as follows.
\begin{equation*}
\begin{aligned}
u^{n+1} &= u^n + \frac{\tau}{2}\left(v^{n} + v^{n+1}\right),\\
v^{n+1} &= v^n - \frac{\tau}{2}A\left(u^{n} + u^{n+1}\right) - \frac{\tau}{2}B\left(v^{n} + v^{n+1}\right) - \frac{\tau}{2}\left(\gamma^n v^{n} +\gamma^{n+1} v^{n+1}\right) + \frac{\tau}{2}\left(f^{n} + f^{n+1}\right).
\end{aligned}
\end{equation*}
By $v^{n+\frac{1}{2}}:=\frac{1}{2}\left(v^{n} + v^{n+1}\right)$, we have
\begin{equation*}
\begin{split}
u^{n+1} &= u^n + \tau v^{n+\frac{1}{2}},\\
v^{n+1} &= v^n - \frac{\tau}{2}A\left(u^{n} + u^{n+1}\right) - \tau Bv^{n+\frac{1}{2}} - \tau\gamma^{n+1} v^{n+\frac{1}{2}} + \frac{\tau}{2} (\gamma^{n+1} - \gamma^{n})v^{n} + \frac{\tau}{2}\left(f^{n} + f^{n+1}\right).
\end{split}
\end{equation*}
In the second equation, we eliminate $ u^{n+1}$ using the first one and obtain
\begin{equation*}
v^{n+1} = v^n - \tau Au^{n} - \frac{\tau^2}{2}Av^{n+\frac{1}{2}} - \tau Bv^{n+\frac{1}{2}} - \tau\gamma^{n+1} v^{n+\frac{1}{2}} + \frac{\tau}{2} (\gamma^{n+1} - \gamma^{n})v^{n}+ \frac{\tau}{2}\left(f^{n} + f^{n+1}\right).
\end{equation*}
Finally, we can write \eqref{CN} to be a half-full-half step formulation
\begin{equation} \label{CN_equiv}
\small
\begin{split}
v^{n+\frac{1}{2}} &= v^n - \frac{\tau}{2} Au^{n} - \frac{\tau^2}{4}Av^{n+\frac{1}{2}} - \frac{\tau}{2} Bv^{n+\frac{1}{2}} - \frac{\tau}{2}\gamma^{n+1} v^{n+\frac{1}{2}} + \frac{\tau}{4} (\gamma^{n+1} - \gamma^{n})v^{n} + \frac{\tau}{4}\left(f^{n} + f^{n+1}\right),\\
u^{n+1} &= u^n + \tau v^{n+\frac{1}{2}},\\
v^{n+1} &= v^{n+\frac{1}{2}}  - \frac{\tau}{2} Au^{n} - \frac{\tau^2}{4}Av^{n+\frac{1}{2}} - \frac{\tau}{2} Bv^{n+\frac{1}{2}} - \frac{\tau}{2}\gamma^{n+1} v^{n+\frac{1}{2}} + \frac{\tau}{4} (\gamma^{n+1} - \gamma^{n})v^{n}+ \frac{\tau}{4}\left(f^{n} + f^{n+1}\right).
\end{split}
\end{equation}
%
%\begin{equation*} \label{}
%\small
%\begin{split}
%v^{n+\frac{1}{2}} &= v^n - \frac{\tau}{2} Au^{n} - \frac{\tau^2}{4}Av^{n+\frac{1}{2}} - \frac{\tau}{2} Bv^{n+\frac{1}{2}} - \frac{\tau}{2}\gamma^{n+\frac{1}{2}} v^{n+\frac{1}{2}} - \frac{\tau}{8} (\gamma^{n+1} - \gamma^{n})(v^{n+1}-v^{n}) + {\color{red}\frac{\tau}{4}\left(f^{n} + f^{n+1}\right)},\\
%u^{n+1} &= u^n + \tau v^{n+\frac{1}{2}},\\
%v^{n+1} &= v^{n+\frac{1}{2}}  - \frac{\tau}{2} Au^{n} - \frac{\tau^2}{4}Av^{n+\frac{1}{2}} - \frac{\tau}{2} Bv^{n+\frac{1}{2}} - \frac{\tau}{2}\gamma^{n+\frac{1}{2}} v^{n+\frac{1}{2}} - \frac{\tau}{8} (\gamma^{n+1} - \gamma^{n})(v^{n+1}-v^{n})+ {\color{red}\frac{\tau}{4}\left(f^{n} + f^{n+1}\right)}.
%\end{split}
%\end{equation*}

\section{vanilla IMEX scheme for Equation~\eqref{eq:evolution}}
\label{sec:IMEX}
If we replace the trapezoidal rule for the nonlinearity in \eqref{CN_equiv} by the left/right rectangular rule, respectively, we have the vanilla IMEX scheme as follows
\begin{align} 
v^{n+\frac{1}{2}} &= v^n - \frac{\tau}{2} Au^{n} - \frac{\tau^2}{4}Av^{n+\frac{1}{2}} - \frac{\tau}{2} Bv^{n+\frac{1}{2}} - \frac{\tau}{2}\gamma^{n+1} v^{n+\frac{1}{2}} + \frac{\tau}{4} (\gamma^{n+1} - \gamma^{n})v^{n}+ \frac{\tau}{2}f^{n}, \label{IMEX_1}\\
u^{n+1} &= u^n + \tau v^{n+\frac{1}{2}}, \label{IMEX_2}\\
v^{n+1} &= v^{n+\frac{1}{2}}  - \frac{\tau}{2} Au^{n} - \frac{\tau^2}{4}Av^{n+\frac{1}{2}} - \frac{\tau}{2} Bv^{n+\frac{1}{2}} - \frac{\tau}{2}\gamma^{n+1} v^{n+\frac{1}{2}} + \frac{\tau}{4} (\gamma^{n+1} - \gamma^{n})v^{n} + \frac{\tau}{2} f^{n+1}, \label{IMEX_3}
\end{align}
which can be regard as the combination of the Crank-Nicolson scheme for the linear part and the leapfrog scheme for the nonlinear part, respectively. This numerical scheme is developed in~\cite{Hochbruck2021imex} for the first time.

The scheme \eqref{CN_equiv} can be also reformulated as
\begin{equation}\label{IMEX_4}
P_{+}\mathbf{u}^{n+1} = P_{-}\mathbf{u}^{n} + \frac{\tau}{2}\left(g^{n} + g^{n+1}\right)+\frac{\tau^2}{4}\left[
\begin{array}{c}
f^{n} - f^{n+1}\\
-(B+\gamma^{n+1})\left(f^n - f^{n+1}\right)
\end{array}
\right].
\end{equation}
Indeed, by subtracting \eqref{IMEX_1} and \eqref{IMEX_3} we have
\begin{equation} \label{IMEX_5}
 v^{n+\frac{1}{2}} = \frac{1}{2}\left( v^n + v^{n+1} \right) + \frac{\tau}{4}\left(f^{n} - f^{n+1}\right).
\end{equation}
Inserting \eqref{IMEX_5} into \eqref{IMEX_2} yields
\begin{equation} \label{IMEX_6}
u^{n+1} = u^n + \frac{\tau}{2}\left( v^n + v^{n+1} \right) + \frac{\tau^2}{4}\left(f^{n} - f^{n+1}\right).
\end{equation}
By combing \eqref{IMEX_1} with \eqref{IMEX_3},  and using \eqref{IMEX_2} and \eqref{IMEX_5} we get
\begin{equation*}
\small
\begin{split}
v^{n+1} &= v^n - \tau Au^{n} - \frac{\tau^2}{2}A v^{n+\frac{1}{2}} - \tau Bv^{n+\frac{1}{2}} - \tau\gamma^{n+1} v^{n+\frac{1}{2}} + \frac{\tau}{2} (\gamma^{n+1} - \gamma^{n})v^{n} + \frac{\tau}{2}\left(f^{n} + f^{n+1}\right)\\
	&= v^n - \tau Au^{n} - \frac{\tau}{2}A\left( u^{n+1} - u^n\right) - \tau Bv^{n+\frac{1}{2}} - \tau\gamma^{n+1} v^{n+\frac{1}{2}} + \frac{\tau}{2} (\gamma^{n+1} - \gamma^{n})v^{n} + \frac{\tau}{2}\left(f^{n} + f^{n+1}\right)\\
	&= v^n - \frac{\tau}{2}A\left( u^{n+1} + u^n\right) - \tau B v^{n+\frac{1}{2}}  - \tau\gamma^{n+1}v^{n+\frac{1}{2}} + \frac{\tau}{2} (\gamma^{n+1} - \gamma^{n})v^{n} + \frac{\tau}{2}\left(f^{n} + f^{n+1}\right)\\
	&= v^n - \frac{\tau}{2}A\left( u^{n+1} + u^n\right) - \frac{\tau}{2}B\left(v^{n} + v^{n+1}\right) - \frac{\tau}{2}\gamma^{n+1}\left(v^{n} + v^{n+1}\right) + \frac{\tau}{2} (\gamma^{n+1} - \gamma^{n})v^{n} + \frac{\tau}{2}\left(f^{n} + f^{n+1}\right) \\
	&\quad - \frac{\tau^2}{4}B\left(f^{n} - f^{n+1}\right) - \frac{\tau^2}{4}\gamma^{n+1}\left(f^{n} - f^{n+1}\right)\\
%	&= v^n - \frac{\tau}{2}A\left( u^{n+1} + u^n\right) - \frac{\tau}{2}B\left(v^{n} + v^{n+1}\right) + \frac{\tau}{2}\left (f^{n} - \gamma^{n}v^{n}  +f^{n+1} - \gamma^{n+1} v^{n+1}\right) \\
%	&\quad - \frac{\tau^2}{4}(B+\gamma^{n+1})\left(f^{n} - f^{n+1}\right).
\end{split}
\end{equation*}
which implies
\begin{equation} \label{IMEX_7}
\small
\begin{split}
v^{n+1} &= v^n - \frac{\tau}{2}A\left( u^{n+1} + u^n\right) - \frac{\tau}{2}B\left(v^{n} + v^{n+1}\right) \\
&\quad+ \frac{\tau}{2}\left (f^{n} - \gamma^{n}v^{n}  +f^{n+1} - \gamma^{n+1} v^{n+1} \right) - \frac{\tau^2}{4}(B+\gamma^{n+1})\left(f^{n} - f^{n+1}\right).
\end{split}
\end{equation}
Therefore, \eqref{IMEX_6} and \eqref{IMEX_7} are two components of \eqref{IMEX_4}.

%
%\begin{align*}
%\small
%P_{+}\mathbf{u}^{n+1} = P_{-}\mathbf{u}^{n} + \frac{\tau}{2}\left(g^{n} + g^{n+1}\right)+\frac{\tau^2}{4}\left[
%\begin{array}{c}
%f^{n} - f^{n+1}\\
%-(B+\gamma^{n+\frac{1}{2}})\left(f^n - f^{n+1}\right)
%\end{array}
%\right]
%+
%\frac{\tau}{4}\left[
%\begin{array}{c}
%0\\
%(\gamma^{n+1} - \gamma^{n})(v^{n+1} - v^{n})
%\end{array}
%\right].
%\end{align*}

\section{Numerical tests for different cases of $\gamma(t)$}
\label{sec:different_gamma}
In the following, we plot the errors of the vanilla IMEX, and the revised IMEX schemes with the fine spatial discretization ($h=0.00697294$) for the problem of Example 1 with another two damping coefficients.
\begin{figure}[!htbp]
\centering
\subfigure[vanilla IMEX scheme]{
%\includegraphics[width=0.4\textwidth]{fig/fig4_coarse.eps}
%\hspace{0.2in}
\includegraphics[width=0.4\textwidth]{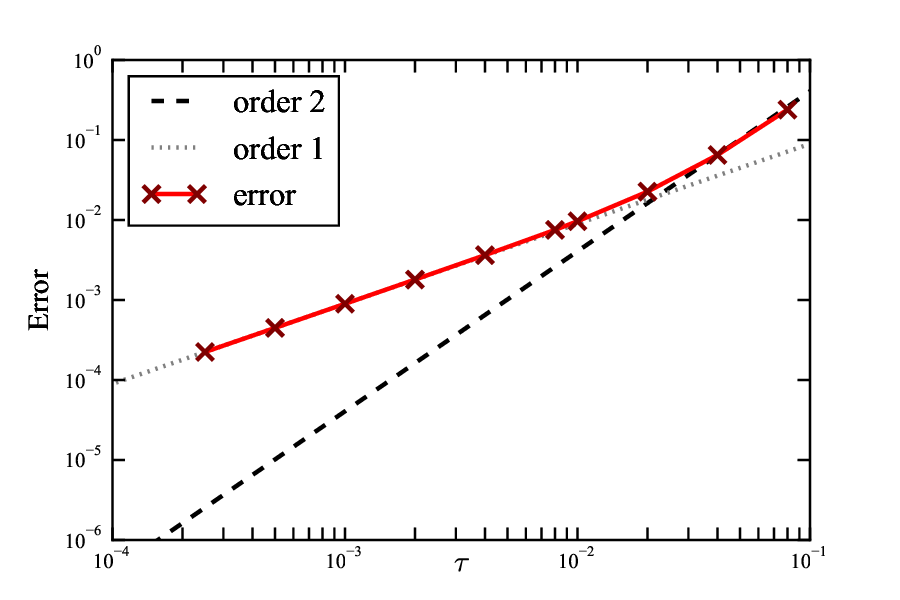}
}
\subfigure[revised IMEX scheme]{
%\includegraphics[width=0.4\textwidth]{fig/fig4_coarse_R.eps}
%\hspace{0.2in}
\includegraphics[width=0.4\textwidth]{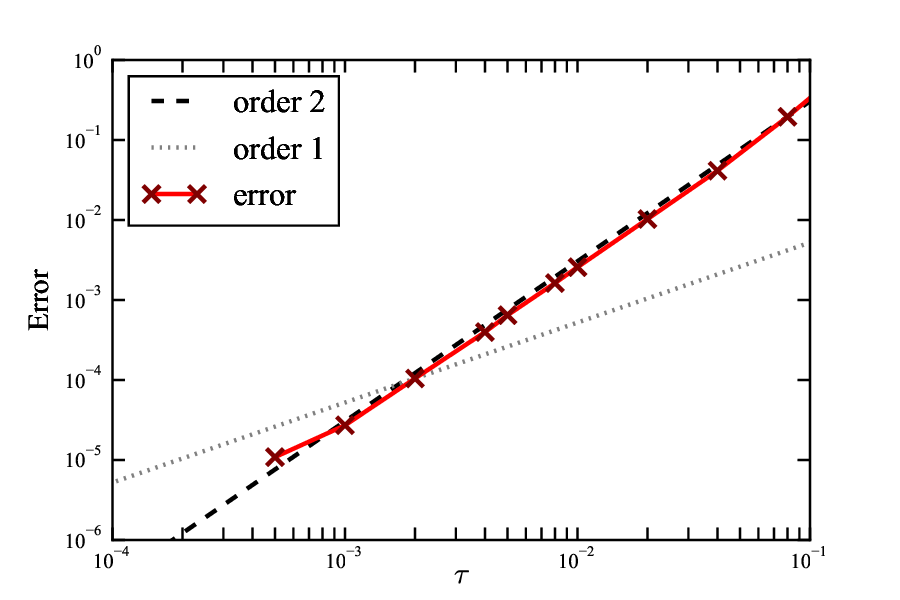}
}
\caption{Error $\mathcal{E}(0.8)$ of the full discretization plots for the problem with the damping coefficient satisfying $r_1 =r_2 =1$ and $\eta = 1$.}
\label{fig_4}
\end{figure} 

\begin{figure}[!htbp]
\centering
\subfigure[vanilla IMEX scheme]{
%\includegraphics[width=0.4\textwidth]{fig/fig5_coarse.eps}
%\hspace{0.2in}
\includegraphics[width=0.4\textwidth]{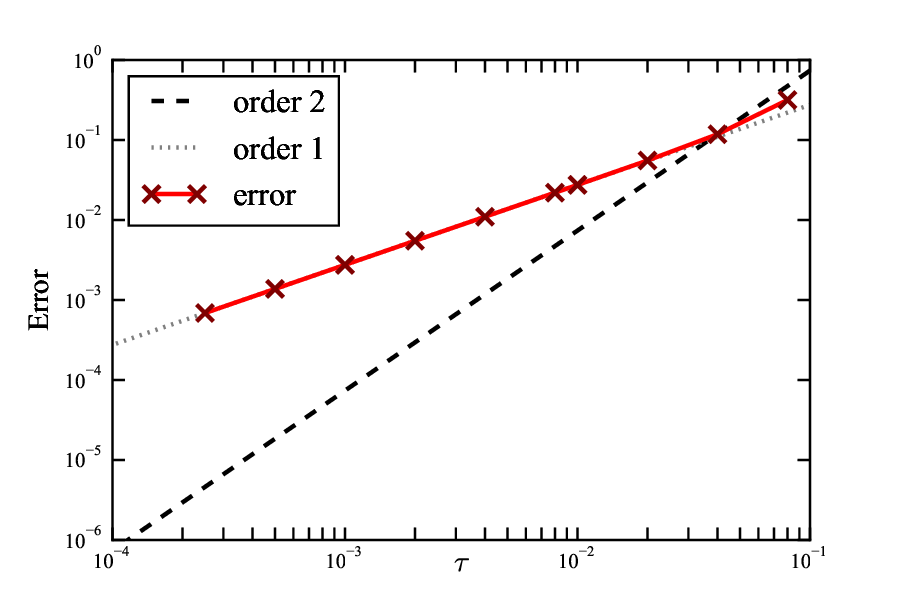}
}
\subfigure[revised IMEX scheme]{
%\includegraphics[width=0.4\textwidth]{fig/fig5_coarse_R.eps}
%\hspace{0.2in}
\includegraphics[width=0.4\textwidth]{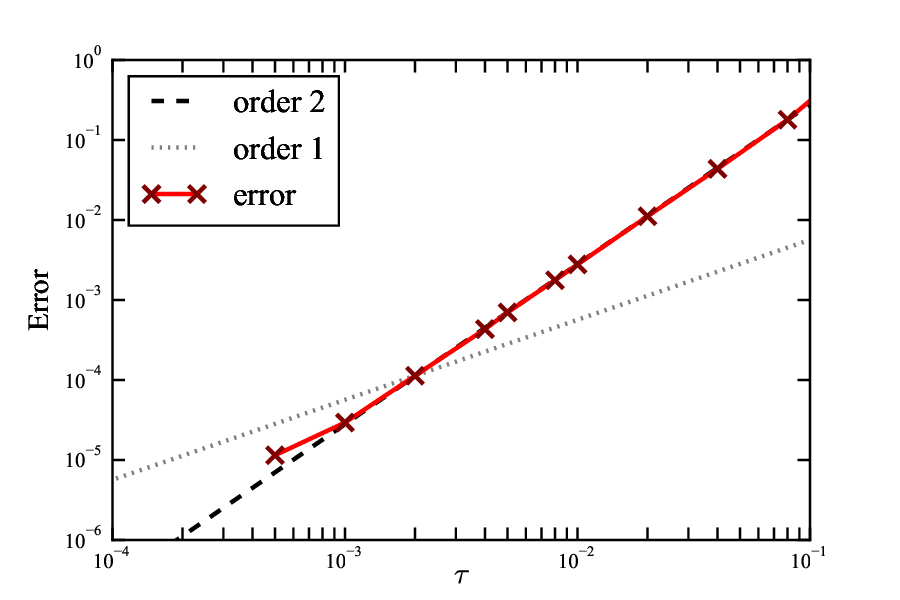}
}
\caption{Error $\mathcal{E}(0.8)$ of the full discretization plots for the problem with the damping coefficient satisfying $r_1 =r_2 = 1$ and $\eta = 2$.}
\label{fig_5}
\end{figure} 
%%%%%%%%%%%%%%%%%%%%%%%%%%%%%%%%%%%%%%%%%%%%%%%%%%%%%%%%%%%%%%%%%%%%%%%%%%%%%%%%
\section*{Acknowledgments}
We thank the anonymous referees very much for the helpful comments.

\section*{Declaration of competing interest}
The authors declare that they have no known competing financial interests or personal relationships that could have appeared to influence the work reported in this paper.

\section*{Data availability}
Data sharing is not applicable to this article as no new data were created or analyzed in this study.

%%%%%%%%%%%%%%%%%%%%%%%%%%%%%%%%%%%%%%%%%%%%%%%%%%%%%%%%%%%%%%%%%%%%%%%%%%%%%%%%

\section*{References}

\bibliographystyle{unsrt}
\bibliography{myref}
\end{document}